\documentclass[journal]{new-aiaa}

\usepackage[utf8]{inputenc}

\usepackage{indentfirst}
\usepackage{algorithm}
\usepackage{algcompatible}

\usepackage{textcomp}
\usepackage{graphicx}
\usepackage{amsmath}
\usepackage[version=4]{mhchem}
\usepackage{siunitx}
\usepackage{longtable,tabularx}
\usepackage{color,soul}
\usepackage{verbatim}
\usepackage[caption=false]{subfig}
\usepackage{float}
\usepackage{booktabs}
\usepackage{multicol}
\usepackage{varioref}
\usepackage{wrapfig}
\usepackage{threeparttable}
\usepackage{dcolumn}
\usepackage{lineno}

\usepackage{amssymb}
\usepackage{mathrsfs}

\graphicspath{{./Fig/}}

\hypersetup{hidelinks}

\newcommand{\de}{\partial}
\newcommand{\di}{{\rm d}}

\newcommand{\vect}[1]{\boldsymbol{#1}}

\title{Robust Bang-Off-Bang Low-Thrust Guidance Using Model Predictive Static Programming}

\author{Yang Wang \footnote{PhD Candidate, Department of Aerospace Science and Technology, Via Masa 32, Milano, 20156, Italy.} and Francesco Topputo\footnote{Assitant Professor, Department of Aerospace Science and Technology, Via Masa 32, Milano, 20156, Italy}}
\affil{Politecnico di Milano, Milano, Italy}

\begin{document}
\footnote[0]{Preliminary results were presented as paper IAC-19,C1,9,3,x49975 at the 70th International Astronautical Congress, Washington D. C., America, 21-25 October 2019. }
\maketitle

Model Predictive Static Programming (MPSP) was always used under the assumption of continuous control, which impedes it for applications with bang-off-bang control directly. In this paper, MPSP is employed for the first time as a guidance scheme for low-thrust transfers with bang-off-bang control where the fuel-optimal trajectory is used as the nominal solution. In our method, dynamical equations in Cartesian coordinates are augmented by the mass costate equation, while the unconstrained velocity costate vector is used as control variable, and is expressed as a combination of Fourier basis functions with corresponding weights. A two-loop MPSP algorithm is designed where the weights and the initial mass costate are updated in the inner loop and continuation is conducted on the outer loop in case of large perturbations. The sensitivity matrix (SM) is recursively calculated using analytical derivatives and SM at switching points is compensated based on calculus of variations. An sample interplanetary CubeSat mission to an asteroid is used as study case to illustrate the effectiveness of the method developed.

\section{Introduction}\label{sec:intro}

In recent decades, highly efficient propulsion systems, such as electric propulsion and solar sails, have made low-thrust engines an alternative to enable ambitious space missions. Extensive work has focused on high-fidelity modeling and open-loop optimal low-thrust trajectory design, solved by direct or indirect methods ~\cite{tang1995,Betts2010,zhang2015,Haberkorn2004}. However, in real-world applications, disturbances such as solar radiation pressure, irregular gravitational fields, outgassing, and unmodeled accelerations, deviate the spacecraft from the nominal trajectory, which requires to update the control profile. The commonly used strategy is to up-link control commands from ground. This requires massive off-line computations and frequent communications between the ground station and the spacecraft. Due to rapid proliferation of space probes, this strategy hardly meets the always-increasing demand for autonomy~\cite{quadrelli2015}.

In literature, a number of guidance laws were proposed for low-thrust orbit transfer problems. Edelbaum~\cite{edelbaum1961propulsion} employed the constant-thrust steering law for quasi-circular orbits. Casalino and Colasurdo~\cite{casalino2007improved} further improved the Edelbaum's method by considering variable specific impulse and thrust magnitude with constant power level. Kluever~\cite{kluever1998simple} designed the simple guidance scheme that blended the individual control law which maximize the time rate of change of a desired orbital element. Petropoulos~\cite{petropoulos2003simple} presented the so-called Q-law for low-thrust orbit transfers, where the proximity quotient Q is served as a candidate Lyapunov function. Hernandez and Akella~\cite{hernandez2014lyapunov,hernandez2016energy} designed Lyapunov control methods for low-thrust orbit transfers Using Levi-Civita and Kustaanheimo-Stiefel transformations. 

Nonlinear optimal control theory (NOCP) is attractive to design the guidance and control law since it can handle constraints while optimizing given performance index. However, the previous mentioned methods are not based on NOCP. Several methods based on NOCP have been designed to track the nominal solution which is computed offline. Neighboring optimal control (NOC) calculates feedback control by optimizing the second-order performance index, which calculates and stores the gain matrix at each time instance off-line and extracts gain matrix using interpolation on-line~\cite{BrysonHo-807}. Pontani et al~\cite{Pontani2015a} proposed variable-time-domain neighboring optimal guidance (VTD-NOG) that avoids numerical difficulty caused by the singularity of gain matrices at terminal time. Zheng~\cite{chen2017a,chen2018b} proposed the backward sweeping algorithm from geometric point of view, for both fixed terminal time and free terminal time of low-thrust transfer problems. Di Lizia et al~\cite{Lizia2008application,Lizia2014high} designed the high-order NOC control law, by using high-order Taylor series automatically achieved by differential algebra around the nominal trajectory. Besides, model predictive control (NMPC) which employs the iterative and finite-horizon optimization strategy has been applied to design the controller. Based on orbit averaging techniques, Gao designed NMPC \cite{gao2008low} to track the mean orbit elements. Huang et al~\cite{huang2012nonlinear} proposed a NMPC strategy using differential transformation based optimization method to track the nominal trajectory. However, the tracking approaches is the overdependence on the nominal profile. For example, the bang-bang thrust sequence is assumed to be unchanged under perturbations when using the NOC method~\cite{chen2017a}. Also, these techniques lack the operational flexibility since the trajectory is restricted to the vicinity of nominal solution. In order to overcome these drawbacks, the algorithms that enable the spacecraft to re-compute the entire nominal trajectory on-line at the beginning of each guidance interval is attractive. There have been several attempts to design the efficient algorithms. Wang~\cite{wang2018minimum} proposed to use convex programming to calculate fuel-optimal spacecraft trajectory. Pesch~\cite{pesch1989a,pesch1989b} used multiple shooting method to re-compute the trajectory for general optimal control problems.

In this work, model Predictive Static Programming (MPSP), an optimal control design technique that combines the philosophy of model predictive control and approximate dynamic programming~\cite{padhi2008}, is designed for low-thrust neighboring control law. The innovation of this method has three aspects \cite{maity2014}. Firstly, it successfully converts a dynamic programming problem to a static programming problem, and thus it requires only a static costate vector for the update of the control profile. Secondly, the symbolic costate vector enables a closed-form solution, which reduces computational load. Thirdly, the sensitivity matrix (SM) which is necessary for the calculation of the static costate vector can be computed recursively. These advantages promote wide applications of the MPSP technique, e.g., terminal guidance~\cite{oza2012}, reentry guidance~ \cite{halbe2013}, and lunar landing guidance~\cite{zhang2016}, etc. Some variations of the MPSP technique have also been proposed to enhance the algorithm performance. For example, the generalized MPSP \cite{maity2014} formulates the problem in continuous-time framework, which does not require any discretization process to begin with. Quasi-Spectral MPSP \cite{mondal2017} expresses the control profile as a weighted sum of basis functions, enabling the method to optimize only a set of coefficients instead of optimizing the control variable at every grid point. However, most works assume continuity of the control profile, which impedes its application for low-thrust transfer missions with bang-off-bang control.

Considering that MPSP technique is an inherent Newton-type method that requires a good initial guess solution \cite{pan2018}, MPSP as a potential neighboring control law for low-thrust transfer problems is investigated in this work. Firstly, the fuel-optimal low-thrust problem is stated in  Cartesian coordinates, where the necessary conditions are formulated based on Pontryagin minimum principle (PMP). The fuel-optimal solution is used as the nominal solution solved by an indirect method. Secondly, inspired by the natural feedback controller given by PMP, the unconstrained costate variable related to the velocity is used as new control variable in MPSP design. In order to ensure the continuity of the switching function at switching points, dynamical equations are augmented by the mass costate equation. Thirdly, SM is recursively calculated using analytical derivative, where SM at switching points is compensated based on calculus of variations. Since SM discontinuity would result in discontinuity of discrete control sequence, the control profile is represented by the combination of Fourier basis functions and corresponding weights, where the weights are initialized based on nominal trajectory using least square method, and updated using Newton's method. Two-loop MPSP algorithm structure is designed for both small and large perturbations, where Newton's method and continuation are implemented in inner and outer loops respectively. The presented MPSP technique is successfully applied to bang-off-bang control for the first time in literature, without resorting to the additional optimization solver. Several numerical simulations are conducted, showing the effectiveness of the proposed method, so enhancing mission flexibility.

This paper is structured as follows: Section 2 states the control problem by using MPSP method. Section 3 depicts the detailed MPSP guidance design. Section 4 presents numerical simulations for a CubeSat mission to an asteroid. Conclusions are given in Section 5.

\section{Problem Statement} \label{sec:statement}
\subsection{Equations of Motion}
This work considers the heliocentric phase of an interplanetary transfer mission. The restricted two-body problem is employed, where the spacecraft subjects to the gravitatonal attraction of the Sun. The spacecraft natural motion consists of Keplerian orbits around the Sun, corresponding to the equation of motion~\cite{battin1999introduction}
\begin{equation} \label{eq1}
    \ddot{\vect r} + \dfrac{\mu}{r^3}\vect r = \vect 0 
\end{equation}
where $\vect r$ is the spacecraft position vector relative to the center of the Sun and $\mu$ is the gravitational parameter. When the low-thrust engine is considered, Eq.\ \eqref{eq1} is modified as 
\begin{equation}
\label{eq:dynamical_eqs}
    \dot{\vect x} = \vect f(\vect x,\vect \alpha,u)
    \Rightarrow
    \begin{pmatrix}
    \dot{\vect r}\\
    \dot{\vect v}\\
    \dot m
    \end{pmatrix}
    =
    \begin{pmatrix}
         \vect v \\
         \vect g(\vect r)+ u \dfrac{T_{\rm max}}{m} \vect \alpha \\
         -u\dfrac{T_{\rm max}}{c}
    \end{pmatrix}
\end{equation}
where $\vect g(\vect r):=-\mu\vect r/r^3$, $\vect r := [x,y,z]^\top \in \mathbb{R}^3$ and $\vect v := [v_x,v_y,v_z]^\top \in \mathbb{R}^3$ are the gravitational vector field, the spacecraft position vector, and its velocity vector, respectively; $m$ is the spacecraft mass, $T_{\rm max}$ is the maximum thrust magnitude, $c = I_{\rm sp} g_0$ is the exhaust velocity ($I_{\rm sp}$ is the engine specific impulse, $g_0$ is the gravitational acceleration at sea level), $u$ is the thrust throttle factor, $\vect \alpha$ is the thrust pointing vector. The state vector is $\vect x = \left[\vect r,\vect v,m \right]\in \mathbb{R}^7$. Both $T_{max}$ and $c$ are assumed constant during flight.

\subsection{Fuel-Optimal Problem}
In this work, fuel-optimal low-thrust trajectory is employed as the reference trajectory. The corresponding performance index is
\begin{equation}
    J = \dfrac{T_{\rm max}}{c}\int_{t_0}^{t_f} u \ \di t
\end{equation}
where $t_0$ and $t_f$ are initial and terminal time instants, both fixed. The initial state is known, i.e., $\vect x(t_0) = \vect x_0$. For the interplanetary mission to the asteroid, the fixed terminal constraint is considered, as 
\begin{equation} \label{eq:terminalconstr}
    \vect x (t_f) = \vect x_f
\end{equation}

The inequality constraint for thrust throttle factor $u$ is
\begin{equation}
    0 \leq u \leq 1
\end{equation}

The Hamiltonian function reads\ \cite{BrysonHo-807}
\begin{equation}
    \label{eq:Hamiltonian}
    H = \dfrac{T_{\rm max}}{c}u + \vect\lambda^T_r \vect v + \vect\lambda^T_v\left[\vect g(\vect r)  + \dfrac{T_{\rm max}}{m}u \vect\alpha\right] - \lambda_m \dfrac{u{T_{\rm max}}}{c}
\end{equation}
where $\vect\lambda = [\vect\lambda_r, \vect\lambda_v, \lambda_m]$ is the costate vector associated with $\vect x$. Dynamical equations of $\vect\lambda$ are
\begin{equation}
    \label{eq:costate_eqs}
    \dot{\vect\lambda} = -\left(\dfrac{\de H}{\de \vect x}\right)^\top 
    \Rightarrow
    \begin{pmatrix}
    \dot{\vect\lambda}_r\\
    \dot{\vect\lambda}_v\\
    \dot \lambda_m
    \end{pmatrix}
    =
    \begin{pmatrix}
         -\vect G^\top \vect\lambda_v \\
         -\vect\lambda_r \\
         u T_{\rm max}/m^2 \vect\lambda^\top_v \vect\alpha
    \end{pmatrix}
\end{equation}
where $\vect G = \de \vect g(\vect r)/\de  \vect r $. Since the final mass is free, there exists
\begin{equation}
\label{eq:lambda_mf}
    \lambda_m (t_f) = 0
\end{equation}

According to PMP, the optimal thrust direction is along the opposite direction of the primer vector $\vect\lambda_v$, as
\begin{equation}
    \label{eq:thrust_dirction}
    \vect\alpha^* = -\dfrac{\vect\lambda_v}{\lambda_v}, \quad  \mbox{if}\quad  \lambda_v \not= 0
\end{equation}
substituting Eq.\ \eqref{eq:thrust_dirction} into Eq.\ \eqref{eq:Hamiltonian} yields
\begin{equation}
\label{eq:hamiltonian_1}
    H =  \vect\lambda^\top_r \vect v + \vect\lambda^\top_v \vect g(\vect r) + \dfrac{u T_{\rm max}}{c}S 
\end{equation}
where the switching function $S$ is defined as
\begin{equation}
\label{eq:switching_func}
    S = -\lambda_v\frac{c}{m} - \lambda_m + 1
\end{equation}

The optimal $u^*$ is governed by $S$ through
\begin{equation} \label{eq:uswitch}
u^* =
\begin{cases} 
0, & \mbox{if} \quad S > 0 \\
1, & \mbox{if} \quad S < 0
\end{cases}
\end{equation}
which is a bang-off-bang control type, forming the thrust sequence.

The fuel-optimal problem is solved by indirect method, which is to find $\vect\lambda_0$ that (together with $\vect x_0$) allows integrating Eqs.\ \eqref{eq:dynamical_eqs} and \eqref{eq:costate_eqs} with the control law in Eqs.\ \eqref{eq:thrust_dirction} and \eqref{eq:uswitch} and verifies the terminal constraints \eqref{eq:terminalconstr} and \eqref{eq:lambda_mf}~\cite{zhang2015}. Singular thrust arcs are not considered here since they have been shown to be non-optimal in general \cite{robbins1965optimality}. Once the optimal $\vect\alpha^*(t)$ and $u^*(t)$ are determined, the spacecraft trajectory can be generated by integrating Eqs.\ \eqref{eq:dynamical_eqs} and \eqref{eq:costate_eqs}. 

\subsection{MPSP Dynamics and Control}

In real world flight, disturbances or new mission requirements need the spacecraft to have the capability update the control sequence automatically. The guidance scheme based on Model predictive static programming (MPSP) is of interest~\cite{padhi2008}. However, the fuel-optimal problem is the optimal control problem with control constraint. MPSP cannot be appied to it directly, since MPSP is originally designed for the unconstrained problem \cite{padhi2008}. In this work, the augmented dynamics and the new control variable are proposed.

Let $\vect x^*(t)$ and $\vect\lambda^*(t)$ denote the reference state and costate profiles, and let $\vect x(t)$ and $\vect\lambda(t)$ be the associated, off-nominal profiles. Let $\Delta\vect\lambda(t)$ be the costate deviation, the two functions~\cite{chen2017a}
\begin{equation}
\label{eq:feedback_contorl}
\begin{cases} 
u(\vect x,\vect\lambda^* + \Delta \vect\lambda) = 1 - \rm{Sgn}\left[ S(\vect x,\vect \lambda^* + \Delta\vect\lambda) \right] \\
\vect\alpha (\vect x,\vect\lambda^* + \Delta \vect\lambda) = -\dfrac{\vect\lambda^*_v + \Delta \vect\lambda_v}{ \lVert \vect\lambda^*_v + \Delta\vect\lambda_v \lVert}, \quad \mbox{if} \ \lVert \vect\lambda^*_v + \Delta\vect\lambda_v \lVert \neq 0
\end{cases}
\end{equation}
define the feedback controller associated with $\vect x$ at time instant $t$, where $\rm{Sgn}$ function is defined as
\begin{equation}
\rm{Sgn}\ (z) =
\begin{cases} 
0,  & \mbox{if} \quad z > 0 \\
1, & \mbox{if} \quad z < 0
\end{cases}
\end{equation}

In this work, the unconstrained costate vector is used as a new control variable for MPSP controller design. This idea also has been utilized in NOC design~\cite{chen2017a} and Lyapunov guidance design \cite{gao2010optimization}. Notice from Eqs.\ \eqref{eq:switching_func} and \eqref{eq:feedback_contorl} that, costate variables which affect $u$ and $\vect\alpha$ are $\vect\lambda_v$ and $\lambda_m$. However, only $\vect\lambda_v$ is seen as the new control variable, based on three facts. Firstly, it can be seen from $\dot\lambda_m$ in Eq.\ \eqref{eq:costate_eqs} that $\vect\lambda_v$ and $\lambda_m$ are dependent, and $\lambda_m$ profile is determined by $\vect\lambda_v$. Secondly, if $\lambda_m$ is also used as a control variable, $\dot{\lambda_m}$ cannot be expressed by Eq.\ \eqref{eq:costate_eqs}. The derivative of the switching function $S$ in Eq.\ \eqref{eq:switching_func}
\begin{equation}
    \dot S = -\dot \lambda_v \dfrac{c}{m} - \lambda_v \dfrac{u T_{max}}{m^2} - \dot \lambda_m
\end{equation}
would hardly be continuous because of the presence of $u$. On the other hand, if $\dot\lambda_m$ is remained as Eq.\ \eqref{eq:costate_eqs}, $\dot S$ becomes simply to be
\begin{equation}
\label{eq:time_dot_S}
    \dot S = -\dot \lambda_v \dfrac{c}{m} 
\end{equation}
which is naturally implicitly dependent on $u$. Thirdly, the second-order differential of $\lambda_m$ w.r.t time is not continuous due to discontinuity of $u$ in $\dot \lambda_m$. Since basis functions are used to approximate the control profile in this work, these may not be appropriate  to efficiently capture the discontinuity \cite{FurfaroMortari-2065}.

Thus, the dynamical equations used for MPSP algorithm design in this work are 
\begin{equation}
\label{eq:arg_eqs}
    \dot{\vect X} 
    = \mathcal{F}(t,\vect X,\vect U)
    \Rightarrow
    \begin{pmatrix}
    \dot{\vect r}\\
    \dot{\vect v}\\
    \dot m \\
    \dot \lambda_m
    \end{pmatrix}
    =
    \begin{pmatrix}
         \vect v \\
         \vect g(\vect r)- u \dfrac{T_{\rm max}}{m\lambda_v} \vect \lambda_v \\
         -u\dfrac{T_{\rm max}}{c}\\
         - u T_{max}/m^2 \lambda_v
    \end{pmatrix}
\end{equation}
where $\vect X = \left[\vect x,\lambda_m \right] \in \mathbb{R}^8$, $\vect U = \vect \lambda_v \in \mathbb{R}^3$, and optimal thrust direction Eq.~\eqref{eq:thrust_dirction} is embedded into Eq.~\eqref{eq:arg_eqs}. The relationships between thrust angles and $\vect U$ are
\begin{equation}
	\begin{cases}
	\vect \alpha =  \arctan(\dfrac{\lambda_{v,2}}{\lambda_{v,1}})\\
	\vect \beta =  \arcsin(\dfrac{\lambda_{v,3}}{\|\vect \lambda_v\|_2})
	\end{cases}
\end{equation}
where $\vect \alpha \in [0,360^\circ]$ is the in-plane angle, $\beta \in [-180^\circ,180^\circ]$ is the out-of-plane angle and $\lambda_{v,i}$ is the $i$th element of $\vect \lambda_v$. Once $\vect X(t)$ and $\vect U(t)$ are determined, the profile of $S$ is decided automatically, which then determines the switching time and thrust sequence. In this work, the task of the MPSP algorithm is to determine a suitable $\Delta \vect U$ and $\Delta \lambda_{m0}$ such that the trajectory of the spacecraft obtained by integrating Eq.\ \eqref{eq:arg_eqs} satisfies the required boundary conditions Eqs.~\eqref{eq:terminalconstr} and \eqref{eq:lambda_mf}, while conducting bang-off-bang control.

\section{MPSP Algorithm Design} \label{sec:method}
\subsection{Sensitive Matrix Calculation}
Different from problems with continuous control profile, dynamical discontinuity happens at switching points. Thus, the trajectory cannot be treated as a whole. In this work, the trajectory is split into multiple segments with switching points located at the boundary of each segment. The time instants at the boundary of each segment are $\{t_0,t_1,\cdots,t_M \}$, where $M$ is the number of total segments, $t_0$ and $t_M = t_f$ are initial and final time, respectively, and $t_k,k=1,2,\cdots,M-1$ are the switching times. Let $\{ t_k^{0_-},t_k^{0_+}, t_k^1,...t_k^{N_k} \}, k=0,1,...,M-1$ denote an evenly-spaced time grid within $[t_k,t_{k+1}]$, where $t_k^{0_-}$ and $t_k^{0_+}$ are the time instants across the impulse. For $k$th segment, $N_k$ is the minimum number of points, such that the time step is just less than prescribed maximum time step $h_{max}$. Note also that $t_k^{N_k} = t_{k+1}^{0_-}$. To ease notation, $N_k$ is denoted $N$. Suppose there is no impulse at initial time, then $t_0^{0_-} = t_0^{0_+} = t_0$.

Consider the $k$th time interval $[t_k,t_{k+1}]$, the discrete system dynamics and the output can be written as
\begin{equation}
    \label{eq:expr_x_k_plus_1}
    \vect X_k^{i+1} =  \vect F_k^i(\vect X_k^i, \vect U_k^i) \quad \vect Y_k^i = \vect O(\vect X_k^i)
\end{equation}
where $\vect Y_k^i$ is the output at the $i$th step, which is the function of $\vect X_k^i$. $\vect F_k^i(\vect X_k^i, \vect U_k^i)$ can be obtained using standard integration formula, such as the Euler method\ \cite{padhi2008}.  High-order integration results in higher accuracy, but larger computational load. In this work, the standard 4th order Runge--Kutta integration is used, see the Appendix for the computation of Eq.\ \eqref{eq:expr_x_k_plus_1}.

The primary objective is to obtain an updated control history $\vect U_k^i$ and initial state $\vect X_0^0$ such that the output $\vect Y$ at terminal time, i.e., $\vect Y^{N}_{M-1}$, reaches to the desired value $\vect Y_{d}$. Writing $\vect Y_{M-1}^N$ about $\vect Y_{d}$ in Taylor series expansion and neglecting high-order terms, the error at the terminal output $\Delta \vect Y_{M-1}^N = \vect Y_{M-1}^N - \vect Y_{d}$ is approximated as
\begin{equation}
    \label{eq:terminal_err}
    \Delta\vect Y_{M-1}^N \cong \di \vect Y_{M-1}^N = \left[ \dfrac{\de \vect Y_{M-1}^N}{\de \vect X_{M-1}^N} \right] \di \vect X_{M-1}^N
\end{equation}

The deduction of $\di \vect X_k^{i+1}$ should consider whether the impulse happens or not. For the interval $[t_k^i,t_k^{i+1}],\ i = 1,2,\cdots,N-1$ where there is no switching point, according to Eq.\ \eqref{eq:expr_x_k_plus_1}, there exists
\begin{equation}
    \label{eq:dx_k_i_plus_1}
    \di \vect X_k^{i+1} = \left[ \dfrac{\de \vect F_k^i}{\de \vect X_k^i} \right] \di \vect X_k^i + \left[ \dfrac{\de \vect F_k^i}{\de \vect U_k^i} \right] \di \vect U_k^i
\end{equation}
where $\di \vect X_k^i$ and $\di \vect U_k^i$ are full differentials of state and control vectors.  

For the interval $[t_k^{0_-},t_k^{1}]$ that contains an impulse, we have
\begin{equation}
\label{eq:dx_k_1}
    \begin{aligned}
    \di \vect X_k^1 & = \dfrac{\de \vect F_k^{0_+}}{\de \vect X^{0_+}_k} \di \vect X^{0_+}_k + \dfrac{\de \vect F_k^{0_+}}{\de \vect U^{0_+}_k} \di\vect U^{0_+}_k \\
    & = \dfrac{\de \vect F_k^{0_+}}{\de \vect X^{0_+}_k}\dfrac{\de \vect X^{0_+}_k}{\de \vect X^{0_-}_k}
    \di\vect X^{0-}_k + \dfrac{\de \vect F_k^{0_+}}{\de \vect X^{0_+}_k}\dfrac{\de \vect X^{0_+}_k}{\de \vect U^{0_-}_k}
    \di\vect U^{0_-}_k + \dfrac{\de \vect F_k^{0_+}}{\de \vect U^{0_+}_k} \di \vect U^{0_+}_k \\
    & = 
    \dfrac{\de \vect F^{0}_k}{\de \vect X^0_k} \di\vect X^0_k + \dfrac{\de \vect F^{0}_k}{\de \vect U^{0}_k} \di\vect U^{0}_k
    \end{aligned}
\end{equation}
where $\di\vect U^{0}_k = \di \vect U^{0_-}_k = \di \vect U^{0_+}_k$ due to thrust angle continuity, $\di \vect X^0_k = \di \vect X^{0_-}_k$, and
\begin{equation}
    \dfrac{\de \vect F^{0}_k}{\de \vect X^0_k} = \dfrac{\de \vect F_k^{0_+}}{\de \vect X^{0_+}_k}\dfrac{\de \vect X^{0_+}_k}{\de \vect X^{0_-}_k}, \quad 
    \dfrac{\de \vect F^{0}_k}{\de \vect U^{0}_k} =  \dfrac{\de \vect F_k^{0_+}}{\de \vect X^{0_+}_k}\dfrac{\de \vect X^{0_+}_k}{\di\vect U^{0_-}_k} + \dfrac{\de \vect F_k^{0_+}}{\de \vect U^{0_+}_k}
\end{equation}

The explressions of ${\de \vect X^{0_+}_k}/{\de \vect X^{0_-}_k}$ and ${\de \vect X^{0_+}_k}/{\de \vect U^{0_-}_k}$ are\ \cite{russell2007primer}
\begin{equation}
    \left[ \dfrac{\de \vect X^{0_+}_k}{\de \vect X^{0_-}_k} \right] = \vect I + \left( \dot{\vect X}^{0_+}_k - \dot{\vect X}^{0_-}_k \right)S_{\vect X}/\dot{S}
\end{equation}

\begin{equation}
    \left[ \dfrac{\de \vect X^{0_+}_k}{\de \vect U^{0_-}_k} \right] = \left( \dot{\vect X}^{0_+}_k - \dot{\vect X}^{0_-}_k \right)S_{\vect U}/\dot{S}
\end{equation}
where $S_{\vect X}$ and $S_{\vect U}$ are row vectors that are the partial derivative of switching function $S$ w.r.t $\vect X$ and $\vect U$ respectively, and $\dot{S}$ is calculated according to Eq.\ \eqref{eq:time_dot_S}. 

Combining Eq.\ \eqref{eq:dx_k_i_plus_1} with Eq.\ \eqref{eq:dx_k_1} yields the uniform form of $\di \vect X_k^{i+1}$ as
\begin{equation}
    \label{eq:dx_unique_form}
    \di \vect X_k^{i+1} = \left[ \dfrac{\de \vect F_k^i}{\de \vect X_k^i} \right] \di \vect X_k^i + \left[ \dfrac{\de \vect F_k^i}{\de \vect U_k^i} \right] \di \vect U_k^i, \quad i=0,1,\cdots,N-1  
\end{equation}

Substituting Eq.\ \eqref{eq:dx_unique_form} into Eq.\ \eqref{eq:terminal_err} yields
\begin{equation}
    \di \vect Y_{M-1}^N = \left[ \dfrac{\de \vect Y_{M-1}^N}{\de \vect X_{M-1}^N} \right] \left \{ \left[ \dfrac{\de \vect F_{M-1}^{N-1}}{\de \vect X_{M-1}^{N-1}} \right] \di\vect X_{M-1}^{N-1} + \left[ \dfrac{\de \vect F_{M-1}^{N-1}}{\de \vect U_{M-1}^{N-1}} \right] \di \vect U_{M-1}^{N-1} \right \}
\end{equation}

Similarly, the state differential at time step $(N-1)$ can be expanded in terms of state and control differentials at time step $(N-2)$. Next, $\di \vect X_{M-1}^{N-2}$ can be expanded in terms of $\di \vect X_{M-1}^{N-3}$ and $\di \vect U_{M-1}^{N-3}$. For $(M-1)$th segment, this process is continued to $\di \vect X_{M-1}^0$. Notice that $\di \vect X_{M-2}^N = \di \vect X_{M-1}^0$, the same process is continued at $(M-2)$th segment. Extending the process until $\vect X_0^0$, one obtains
\begin{equation}
    \label{eqs:dyn_expansion}
    \begin{aligned}
    \di \vect Y^N_{M-1} &= \vect A \di\vect X^0_0 + \vect B_0^0\di\vect U_0^0 + \vect B_0^1\di\vect U_0^1  + \vect B_{M-1}^{N-1}\di\vect U_{M-1}^{N-1} \\ 
    & = \vect A\di \vect X_0^0 + \sum_{k=0}^{M-1} \sum_{i=0}^{N-1} \vect B_k^i\di \vect U_k^i
    \end{aligned}
\end{equation}
where the compact form of coefficients $\vect A$ and $\vect B_k^i$ in Eq.\ \eqref{eqs:dyn_expansion} are
\begin{equation}
    \label{eq:mpap_coeffs}
    \begin{array}{cc}
         \vect A =  \left[ \dfrac{\de \vect Y_{M-1}^N}{\de \vect X_{M-1}^N} \right] \prod\limits_{k=0}^{M-1} \prod\limits_{i=0}^{N-1} \left[ \dfrac{\de \vect F_k^i}{\de \vect X_k^i} \right] \\
         \vect B^i_k  = \left[ \dfrac{\de \vect Y_{M-1}^N}{\de \vect X_{M-1}^N} \right] \left\{\prod\limits_{p=M-1}^{k+1} \prod\limits_{q=N-1}^{1}\left[ \dfrac{\de \vect F_p^q}{\de \vect X_p^q} \right]\right\} 
        \left\{\prod\limits_{q=N-1}^{i+1}\left[ \dfrac{\de \vect F_k^q}{\de \vect X_k^q} \right]\right\}\dfrac{\de \vect F_k^i}{\de \vect U_k^i}
    \end{array}
\end{equation}

The presented MPSP is desirable because the computation of the sensitivity matrix $\vect B_k^i$ can be reduced to an iterative calculation. Define
\begin{equation}
    \vect B_{M-1,0}^{N} = \left[ \dfrac{\de \vect Y_{M-1}^N}{\de \vect X_{M-1}^N} \right]
\end{equation}
there exists
\begin{equation}
    \vect B^i_{k,0} = \vect B^{i+1}_{k,0} \left[ \dfrac{\de \vect F^{i+1}_k}{\de \vect X^{i+1}_k} \right], \quad \vect B^i_k = \vect B^i_{k,0} \left[ \dfrac{\de \vect F_k^i}{\de \vect U_k^i} \right]
\end{equation}

\subsection{Control Representation and Update}
For applications with continuous control profile, discrete control sequence works due to the continuity of SM. However, $\vect B_k^i$ is discontinuous before and after impulse, resulting in the discontinuity of the thrust angle sequence if the control profile is discretized, which is meaningless from the physical point of view. In this work, the control profile is expressed by basis functions. The advantages lies in two facets. Firstly, the continuity of the thrust angle profile can be ensured automatically due to the continuity of basis functions. Secondly, the time derivative of the switching function $S$ can be calculated analytically. The control is expressed as
\begin{equation}
\label{eq:CF_1}
\vect U_k^i(\eta) = \vect P_k^i(\eta) \vect\epsilon
\end{equation}
where $\epsilon$ is the weight to the basis functions, and	
\begin{equation}
\vect P(\eta) = 
\begin{bmatrix}
\vect h^\top(\eta) & & &\\
& \vect h^\top(\eta) & &\\
& & & \vect h^\top(\eta)
\end{bmatrix}
\end{equation}
$\vect h(\eta)$ is the collection of different orders of basis functions. The $\eta$ range is determined by the basis functions chosen. The linear projection of $\eta$ w.r.t time $t$ is used as
\begin{equation}
\eta = \dfrac{\eta_f-\eta_0}{t_f-t_0}\left(t-t_0\right) + \eta_0
\end{equation}

From Eq.~\eqref{eq:CF_1}, the update of the control profile is achieved by updating $\vect \epsilon$. The differential of Eq.\ \eqref{eq:CF_1} is
\begin{equation}
\label{eqs:diff_control_profile}
\di\vect U_k^i = \vect P_k^i(\eta) \di\vect\epsilon
\end{equation}

Substituting Eq.\ \eqref{eqs:diff_control_profile} into Eq.\ \eqref{eqs:dyn_expansion} yields
\begin{equation}
\label{eqs:dYn_expression_new}
        \di\vect Y_{M-1}^N =  \vect A\di\vect X_0^0 + \sum_{k=0}^{M-1} \sum_{i=0}^{N-1} \vect B_k^i\vect P_k^i\di\vect\epsilon
\end{equation}

Since the control used is not directly linked to the thrust throttle factor $u$, the solution in the neighborhood of the reference solution is preferred. The performance index is set to
\begin{equation}
    J = \dfrac{1}{2} \di\vect\epsilon^\top \vect R_{\epsilon} \di\vect\epsilon + \dfrac{1}{2} (\di\vect X_0^0)^\top \vect R_{0} \di\vect X_0^0
\end{equation}

Denote $\vect B_v = \sum_{k=0}^{M-1} \sum_{i=0}^{N-1} \vect B_k^i\vect P_k^i$, the augmented performance index reads
\begin{equation}
      \hat{J} = \dfrac{1}{2} \di\vect\epsilon^\top \vect R_{\epsilon} \di\vect\epsilon + \dfrac{1}{2} (\di \vect X_0^0)^\top \vect R_{0} \di\vect X_0^0  + \vect p^\top\left( \di \vect Y_{M-1}^N - \vect A\di\vect X_0^0 - \vect B_v \di\vect\epsilon \right)
\end{equation}
where $\vect p$ is the associated static costate vector. The optimal conditions read
\begin{equation}
\label{eqs:optimal_cond}
    \begin{array}{cc}
    \left(\dfrac{\di\hat{J}}{\di\left(\di \vect \epsilon\right)}\right)^\top &= \vect R_{\epsilon} \di\vect\epsilon - \vect B_v^\top \vect p = \vect 0 \\
    \left(\dfrac{\di\hat{J}}{\di(\di\vect X_0^0)}\right)^\top &= \vect R_{0} \di \vect X_0^0 - \vect A^\top \vect p = \vect 0\\
    \end{array}
\end{equation}

Substituting Eq.\ \eqref{eqs:optimal_cond} into Eq.\ \eqref{eqs:dYn_expression_new} yields
\begin{equation}
\label{eqs:p_expression}
\vect p = \left(\vect A \vect R_{0}^{-1} \vect A^\top + \vect B_v \vect R_{\epsilon}^{-1} \vect B_v^\top  \right)^{-1}\di\vect Y_{M-1}^N
\end{equation}
and substituting Eq.\ \eqref{eqs:p_expression} into Eq .\ \eqref{eqs:optimal_cond} yields the Newton direction $\di\vect\epsilon$ and $\di\vect X_0^0$, as
\begin{equation}
	\begin{array}{cc}
	\di\vect\epsilon &= \vect R_{\epsilon}^{-1}\vect B_v^\top \vect p\\
	\di \vect X_0^0 & = \vect R_{0}^{-1}\vect A^\top \vect p
	\end{array}
\end{equation}

Since the spacecraft initial position, velocity and mass are known and fixed, $\di\lambda_{m0}$ is to be solved. $\vect \epsilon$ and $\lambda_{m0}$ are updated at $j$th iteration as
\begin{equation}
\begin{array}{cc}
    \vect\epsilon_{j+1} &= \vect\epsilon_{j} - \kappa \di \vect\epsilon_{j}\\
    \lambda_{m0,j+1} &= \lambda_{m0,j} -  \kappa \di\lambda_{m0,j}
\end{array}
\end{equation}
where $\kappa$ is the Newton step length. Then the updated $\vect U(t)$ and $\vect X(t)$ are determined through calculating Eq.\ \eqref{eq:CF_1} and integrating Eq.\ \eqref{eq:arg_eqs}.

However, particular attention needs to be paid for Newton step length $\kappa$. If $\kappa$ is too large, the thrust sequence may be apparently changed, which amplifies undesired terminal error and further deteriorates algorithm performance. On the other hand, the unconstrained and changeable thrust sequence contributes toward enlarging the convergence domain. Therefore, rigorous strategy for $\kappa$ selection should be designed to maintain the algorithm stability and simultaneously enlarge the convergence domain. In this work, the thrust sequence at each iteration is checked and restricted. Denote $N_{seg,i}$ and $N_{seg,ref}$ as the sum of thrust segments and coast segments for the trajectory at $i$th iteration and reference trajectory, respectively. For $(i+1)$th iteration, the restriction on $N_{seg,i+1}$ is implemented such that $\|N_{seg,i+1}-N_{seg,ref}\| \leq N_{\rm seg,tol}$, where $N_{\rm seg,tol}$ is the tolerance for the varied segments. Otherwise, $\kappa$ is reduced. 

\subsection{Nominal Solution Generation}
The fixed-time fuel-optimal open-loop problem can be formulated into a two-point boundary value problem (TPBVP) and it requires to search a zero of the shooting function associated with TPBVP \cite{zhang2015}. In this work, the method that combines analytic derivatives, switching detection technique and numerical continuation is applied to find the fuel-optimal low-thrust trajectory\cite{zhang2015}, which is used as nominal solution. There is no need to assign control structure a prior, and it is also useful in cases where very low-thrust accelerations are used in highly nonlinear vector fields.

The nominal discrete control sequence at $k$th evenly time grid is denoted as $\vect U_{k,ref}$. Collecting all discrete points yields the following linear algebraic equation
\begin{equation}
    \label{eq:Algebric_eqs}
    \vect A \vect\epsilon_0 = \vect B
\end{equation}
where 
\begin{equation}
    \vect A = \left[ \vect P_1, \vect P_2, \cdots, \vect P_N \right]
\end{equation}
\begin{equation}
    \vect B = [\vect U_{1,ref}, \vect U_{2,ref}, \cdots, \vect U_{N,ref}]
\end{equation}

The least-square solution is used as nominal solution, as
\begin{equation}
    \vect\epsilon_0 = (\vect A^\top \vect A)^{-1}\vect A^\top \vect b
\end{equation}

\subsection{Implementation}
The MPSP algorithm for bang-off-bang low-thrust transfers requires to detect the switching time accurately, which is based on two reasons. Firstly, if the switching time is not detected, the integration error will be accumulated around the switching points, which deteriorates the performance of the Newton's method. Secondly, SM is discontinuous as shown in Eqs.\ \eqref{eq:dx_k_i_plus_1} and \eqref{eq:dx_k_1}, which requires to detect switching time for accurate calculation of SM. The switching detection technique is embedded into four-order Runge-Kutta fixed step trajectory integral scheme. The detection is active as soon as the switching function Eq.\ \eqref{eq:switching_func} traverses zero at time interval $[t_k,t_{k+1}]$. The bisection method is used to find the switching time $t_{sw} \in [t_k,t_{k+1}]$ such that the absolute value of switching function is within the tolerance $10^{-12}$.

Based on the techniques proposed above, a two-loop MPSP algorithm is designed consisting of inner-loop and outer-loop parts. The inner-loop algorithm is illustrated in Algorithm \ref{algo:inner-MPSP}, which is to update $\vect \varepsilon$ and $\lambda_{m,0}$ for the input boundary conditions using Newton's method. In inner-loop algorithm, the Newton's method is implemented only when the L2-norm of terminal error at $j$th iteration is less than a maximum error tolerance $\Delta_{\rm max}$, i.e., $\|\vect Y_{M-1}^N\|_{2,j} \leq \Delta_{\rm max}$. $\Delta_{\rm max}$ is the conservative value indicating the failure of the iteration or the encounter of large perturbations. $Sign$ is used to label the success ($Sign = 1$) or failure ($Sign = 0$) of the inner-loop part. The failure occurs when the terminal error exceeds the tolerance or the step length is small enough. The outer-loop MPSP algorithm is shown in Algorithm \ref{algo:outer-MPSP}, which is triggered when $Sign = 0$ is returned. In this case, the continuation from nominal conditions to perturbed conditions is conducted. Denote $\vect C_{\rm ref}$ as reference boundary conditions, $\vect C_{\rm per}$ as perturbed conditions and $\tau$ as continuation parameter. Starting from $\tau = 0$ which corresponds to $\vect C_{\rm ref}$, continuation proceeds until $\tau = 1$ which corresponds to $\vect C_{\rm per}$. At each step, the inner-loop MPSP algorithm is applied to find the solution corresponding to the conditions
\begin{equation}
\vect C_{\tau} = (1-\tau)\vect C_{\rm ref} + \tau \vect C_{\rm per}
\end{equation} 

Beside, the $N_{\rm seg,tol}$ used in the inner-loop algorithm is initially set to be $0$ at the outer-loop algorithm. Thus, the MPSP algorithm tries to find the solution with the same thrust sequence as the nominal solution first. The value of $N_{\rm seg,tol}$ increases once the inner-loop MPSP algorithm fails.

\begin{algorithm}[htp]
	\caption{Inner-Loop MPSP Algorithm}
	\label{algo:inner-MPSP}
	\begin{algorithmic}[1]
			\REQUIRE{Boundary conditions $\vect C_{\tau}$, weight of control sequence $\vect \varepsilon$, $\Delta_{\rm max}$, $N_{\rm seg,tol}$, reference solution}.
			\ENSURE{The updated $\vect \varepsilon$, $\lambda_{m0}$ and label $Sign$}
			\STATE{Integrate the trajectory using initial condition $\vect X_0$ and $\vect \varepsilon$.}
			\STATE{Set $Sign = 1$, $i = 0$.}
			\WHILE{The terminal error does not satisfy requirement}
		    \IF{$\| \vect Y_{M-1}^N \|_{2,i} \geq \Delta_{\rm max}$}
		    \STATE{$Sign = 0$. Return.}
			\ENDIF
			
			\STATE{ Calculate $\vect A$ and $\vect B_v$ matrix in Eq.\ \eqref{eq:mpap_coeffs} and the static costate vector $\vect p$ in Eq.\ \eqref{eqs:p_expression}. Calculate Newton direction $\di \vect \varepsilon_j$ and $\di \lambda_{m0,j}$. Set initial $\kappa = 1$.}
			
			\WHILE{1}
			\STATE{$\vect \varepsilon_{\rm old} := \vect \varepsilon$ and $\lambda_{m,0,{\rm old}} := \lambda_{m,0}$.}
			\STATE{Update $\vect \varepsilon$ and $\lambda_{m,0}$. Integrate dynamical equations to obtain the trajectory. Get $N_{{\rm seg},i+1}$. $i:=i+1$}
			\IF{$\|N_{{\rm seg},i+1}-N_{\rm seg,ref}\| \leq N_{\rm seg,tol}$}
			\STATE{Save the updated control. Break.}
			\ELSE
			\STATE{$\kappa := \kappa/2$. $\vect \varepsilon := \vect \varepsilon_{\rm old}$. $\lambda_{m,0} := \lambda_{m,0,{\rm old}}$.}
			\IF{$\kappa \leq 1/2^5$}
			\STATE{Set $Sign = 0$. Return.}
			\ENDIF
			\ENDIF
			\ENDWHILE
			\ENDWHILE

	\end{algorithmic}
\end{algorithm}

\begin{algorithm}[htp]
	\caption{Outer-Loop MPSP Algorithm}
	\label{algo:outer-MPSP}
	\begin{algorithmic}[1]
		
		\STATE{Solving the fuel-optimal low-thrust transfer problem using indirect method~\cite{zhang2015}. The initial weight $\vect\epsilon_0$ is calculated using least square method.}
		
		\STATE{Calculate perturbed condition or thruster parameters.}
		
		\STATE{Set error tolerance $\Delta_{max}$, continuation parameter variation $\delta \tau$ and $N_{\rm seg,tol}=0$.}
		
		\STATE{Apply inner-loop MPSP algorithm. Denote the returned label of success as $Sign1$.}
		
		\IF{$Sign1 = 0$}
		\WHILE{1}
		
			\STATE{Set $\tau = \delta \tau$, $\tau_{\rm old} = 0$.}
		
			\WHILE{$\delta \tau > 0$}
		
				\STATE{Calculate the perturbed condition for current $\tau$.}
				\STATE{Apply inner-loop MPSP algorithm. Denote the returned label of success as $Sign2$.}
				\IF{$Sign2 = 0$}
				\STATE{$\delta \tau := \delta \tau/2$.}
				\ELSE
				\STATE{Save the solution as new initial guess solution for the next iteration.}
				\STATE{Set $\delta \tau := \min(1-\tau, 2\delta \tau)$, $\tau_{\rm old} = \tau$}
				\ENDIF
				\STATE{$\tau = \tau_{\rm old} + \delta \tau$.}
				\IF{$\delta \tau \leq 0.01$ and $\tau_{\rm old} \neq 1$}
				\STATE{Break}
				\ENDIF
			\ENDWHILE
			
			\IF{$\delta \tau = 0$}
				\STATE{break}
			\ELSE
				\STATE{$N_{\rm seg,tol} : = N_{\rm seg,tol} + 2$}
				\IF{$N_{seg,ref} - N_{\rm seg,tol} < 0$}
				\STATE{Print failure information; break}
				\ENDIF
			\ENDIF

		\ENDWHILE
		\ELSE
			\STATE{Save the solution.}
		\ENDIF
	\end{algorithmic}
\end{algorithm}

\section{Numerical Simulations} \label{sec:simulations}

\subsection{Fuel-optimal Trajectory}
\begin{figure}
	\centering
	\scalebox{0.7}{\includegraphics{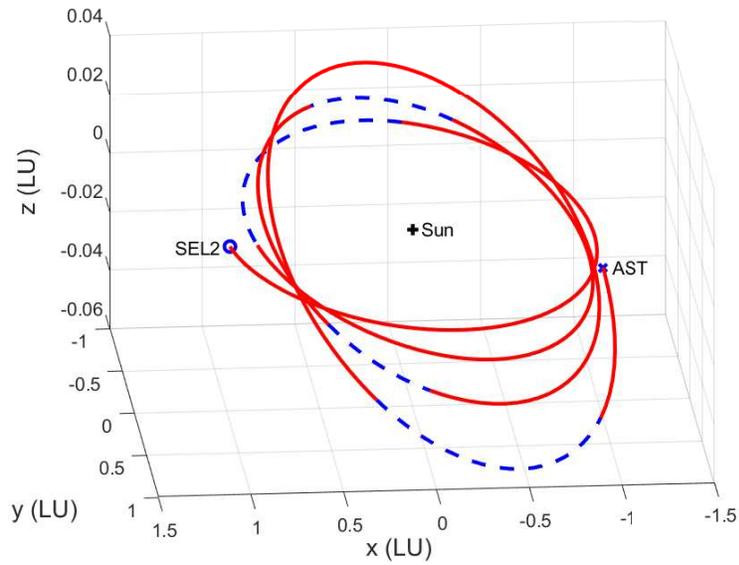}}
	\caption{The fuel-optimal trajectory for boundary conditions in Tab.~\ref{table:boundary_conditions}, where 'SEL2' denotes Sun-Earth L2 starting point and 'AST' denotes the asteroid position.}
	\label{fig:open_tra_ast}
\end{figure}

\begin{figure}
	\centering
	\scalebox{0.7}{\includegraphics{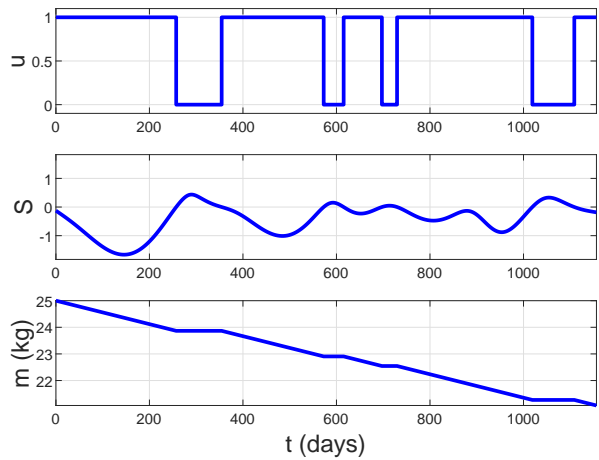}}
	\caption{The variations of thrust throttle $u$, switching function $S$ and mass $m$ w.r.t time corresponding to the fuel-optimal trajectory in Fig.~\ref{fig:open_tra_ast}.}
	\label{fig:open_usm_ast}
\end{figure}

An interplanetary Cubesats mission to the asteroid $99942$ Apophis is considered. The related physical constants are listed in Tab.\ \ref{table:physical_constant}. The initial mass $m_0$, maximum thrust magnitude $T_{\rm max}$ and specific impulsive $I_{\rm sp}$ are set to be $25 kg$, $1.5\times 10^{-3} N$ and $3000 s$, respectively. The departure epoch is October 1st, 2020, the arrival epoch is December 1st, 2023, and the departure position is Sun-Earth L2. The boundary conditions of the spacecraft at the initial and rendezvous time are listed in Tab.~\ref{table:boundary_conditions}. The fuel-optimal solution is employed as the nominal solution which is solved by using indirect method \cite{zhang2015}. The optimal transfer orbit is shown in Fig.\ \ref{fig:open_tra_ast}, where the red line denotes the thrust segment, i.e., $u=1$, while the blue dash line denotes the coast segment, i.e., $u=0$. The variations of thrust throttle $u$, switching function $S$ and the mass $m$ w.r.t time are shown in Fig.\ \ref{fig:open_usm_ast}. The optimal trajectory consists of five thrust segments and four coast segments, and the final mass of the spacecraft is $21.062 kg$.


In the following numerical simulations, Fourier basis polynomials with maximum $15$th order are used to approximate the control sequence. The largest time step for each segment is set to be $h_{max} = 0.0005 \times t_f$. The convergence conditions to terminate the MPSP algorithm is such that the terminal position error $\|\Delta \vect r_f\|_2 \leq 500 km$, the terminal velocity error $\|\Delta \vect v_f\|_2 \leq 0.1 km/s$ and $|\lambda_{mf}| \leq 10^{-6}$. For other parameter setting, $\tau_0 = 0.5$, $\delta \tau = 0.5$ and $\Delta_{max} = 1$. All simulations are conducted under Intel Core i7-9750H, CPU@2.6GHz, Windows 10 system within MATLAB environment. Considering the paper length, the simulations with large number of random perturbations are not reported.

\begin{table}
	\centering
	\caption{Physical constants.}
	\begin{tabular}{lccc}
		\toprule
		\toprule
		Physical constant & Values \\
		\midrule
		Mass parameter $\mu$ & $1.327124 \times 10^{11} \ km^3/s^2$ \\
		Gravitational field, $g_0$ & $9.80655 \ m/s^2$ \\
		Length unit, LU& $1.495979 \times 10^8 \ km$\\
		Time unit, TU & $5.022643 \times 10^6 \ s$\\
		Velocity unit, VU & $29.784692 \ km/s$\\
		Mass unit, MU & $25 \ kg$ \\
		\bottomrule
	\end{tabular}
	\label{table:physical_constant}
\end{table}

\begin{table}
	\centering
	\caption{Boundary Conditions.}
	\begin{tabular}{lccc}
		\toprule
		\toprule
		Boundary Condition & Values \\
		\midrule
		Initial position vector (LU)  & $\vect r_0 = [1.001367,0.140622,--6.594513\times 10^{-6}]^\top$\\
		Initial velocity vector (VU)  & $\vect v_0 = [-0.155386,0.986258,-4.827818\times 10^{-5}]^\top$ \\
		Terminal position vector (LU) & $\vect r_f = [-1.044138,-0.122918,-0.018183]^\top$\\
		Terminal velocity vector (VU) & $\vect v_f = [0.222668,-0.875235,0.051944]^\top$\\
		\bottomrule
	\end{tabular}
	\label{table:boundary_conditions}
\end{table}

\subsection{Perturbations on Initial Conditions}


In this section, the proposed MPSP algorithm is tested by assuming the perturbations on the initial conditions. Different scales of perturbation magnitude for the single base perturbation are used for analysis. The random base perturbation is set to be as $\delta \vect x_0 = [1.6712,-1.0659,-4.1460,1.0876,2.3763,4.6091] \times 10^{-2}$. The $9$ cases of different scales of perturbation magnitude are simulated, and the corresponding perturbed initial condition is $\delta \hat{\vect x}_0 = \vect x_0 + \kappa \delta \vect x_0$ where $\kappa \in [-3,-2.5,-2,-1,1,2,2.5,3,3.5]$. The perturbation direction of the former $4$ cases are opposite to that of the later $5$ cases. The simulation results are summarized in Tab.~\ref{table:ini_per} which gives the terminal errors, the Newton's iteration in the inner-loop MPSP algorithm and the percentage of fuel increase w.r.t the corresponding optimal solutions. It can be observed that, the algorithm works successfully since the terminal errors are all within the tolerance. 

\begin{figure}
	\centering
	\subfloat[]{\includegraphics[width=0.5\textwidth]{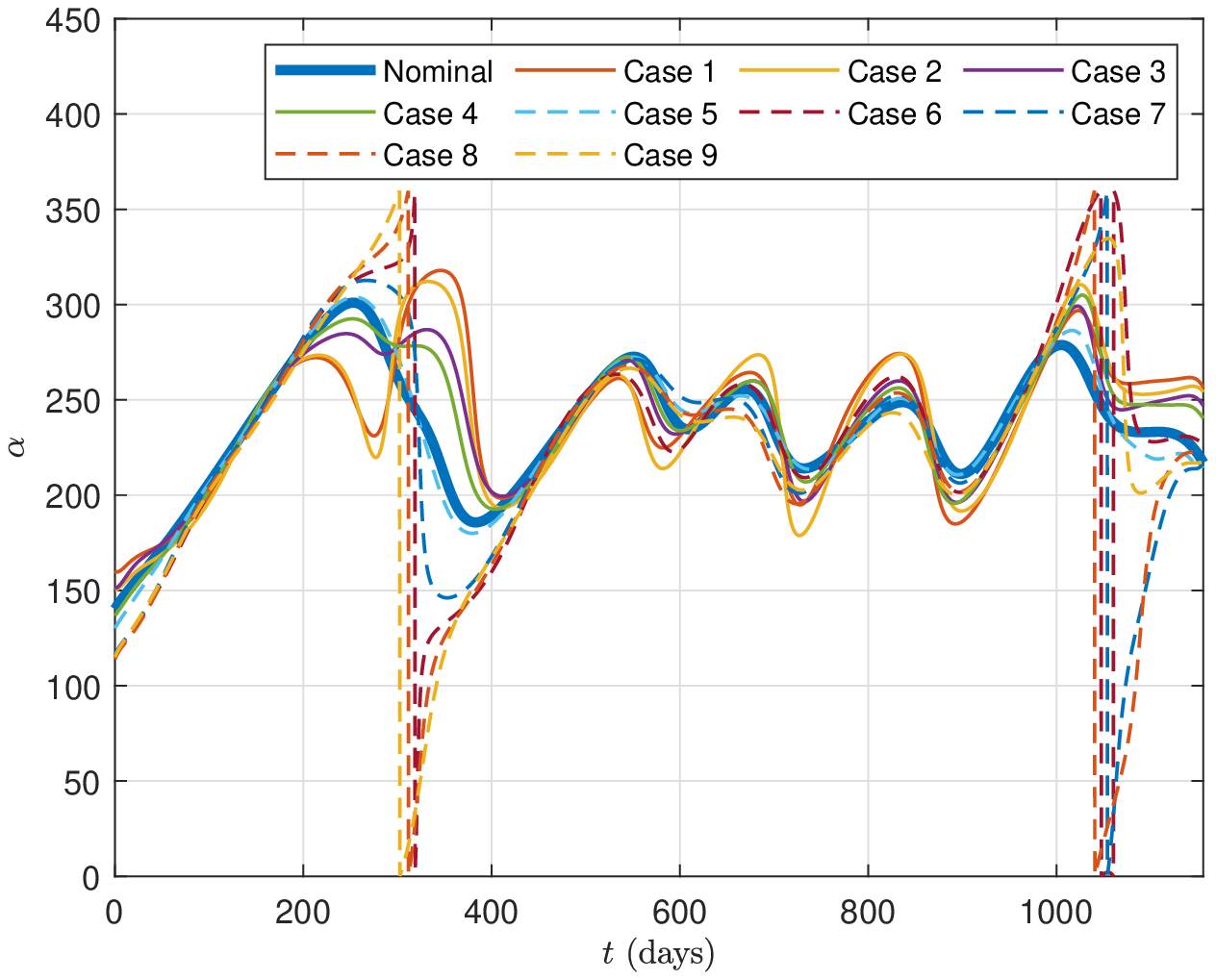}\label{fig:Ini_per_alpha_9_cases}}
	\hfill
	\subfloat[]{\includegraphics[width=0.5\textwidth]{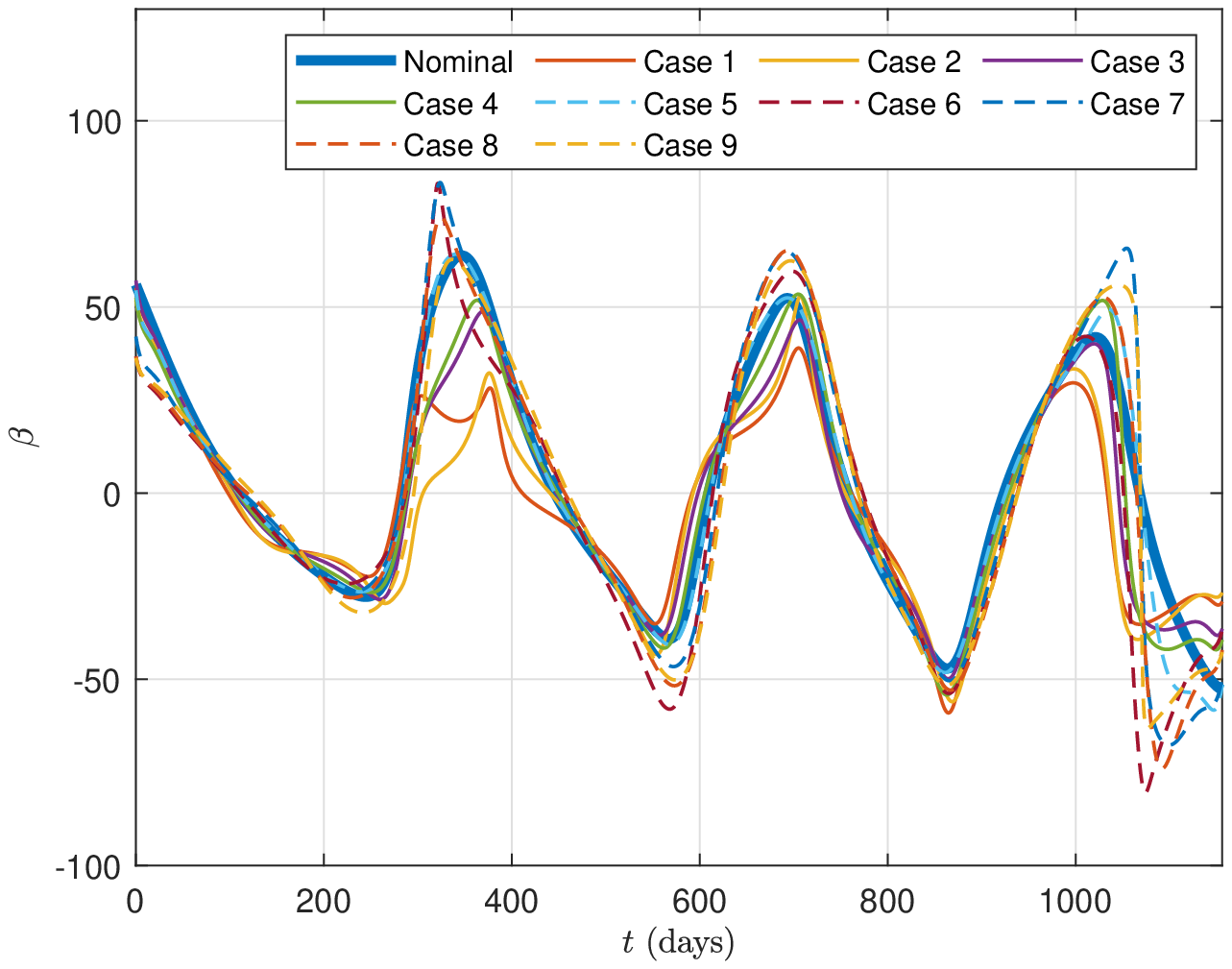}\label{fig:Ini_per_beta_9_cases}}
	
	\caption[Simulation results]{The in-plane angle $\alpha$ and out-of-plane angle $\beta$ of the converged MPSP trajectories for variations of initial conditions in Tab.~\ref{table:ini_per}.}
	\label{fig:Ini_per_9_cases}
\end{figure}

\begin{figure}
	\centering
	{\includegraphics[width=0.7\textwidth]{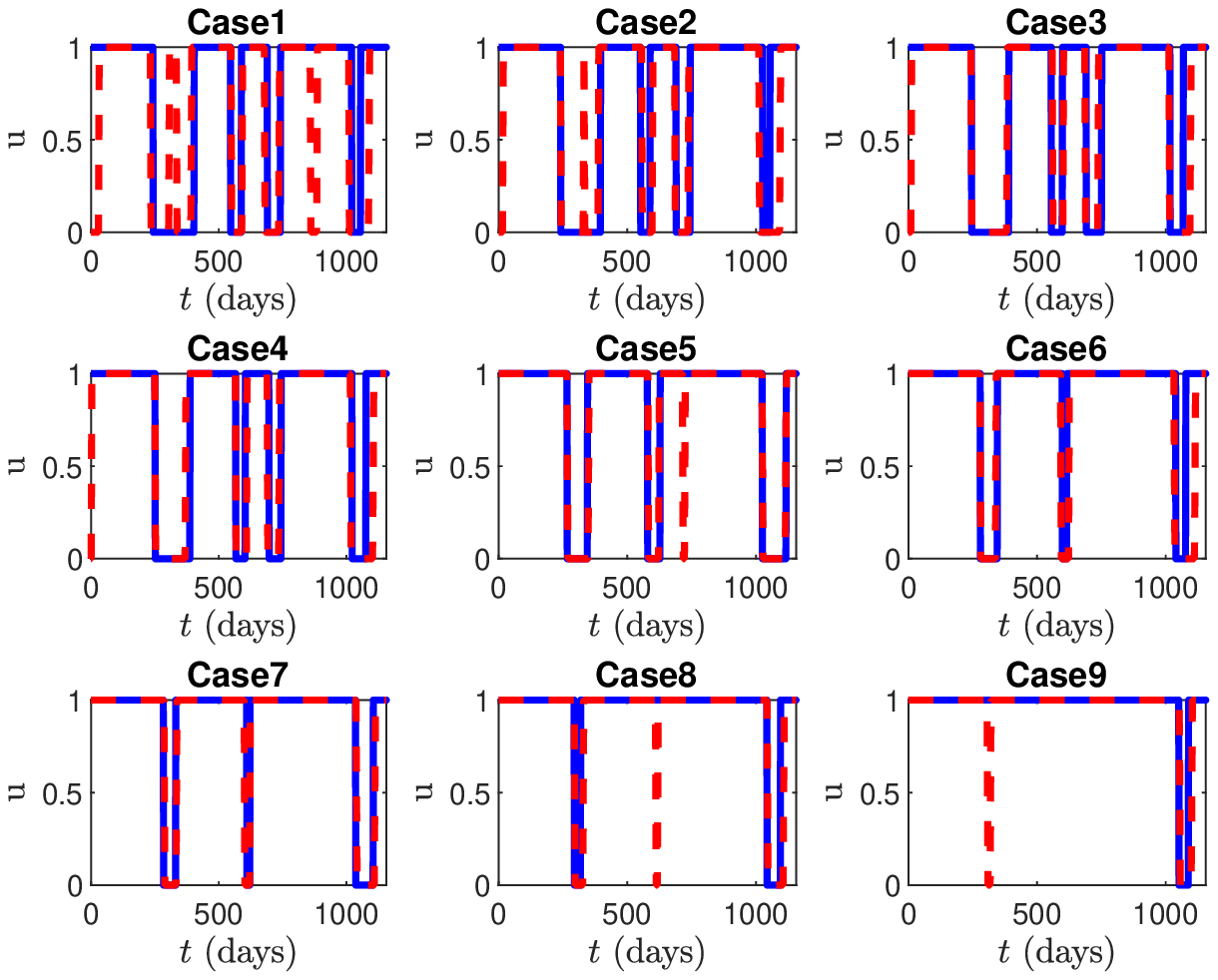}}
	\caption{The comparison of optimal thrust sequences and thrust sequences of MPSP solutions for cases in Tab.~\ref{table:ini_per}. The red dash line: optimal thrust sequence; Blue line: thrust sequences of MPSP solutions.}
	\label{fig:Ini_per_u_9_cases}
\end{figure}

\begin{figure}
	\centering
	{\includegraphics[width=0.7\textwidth]{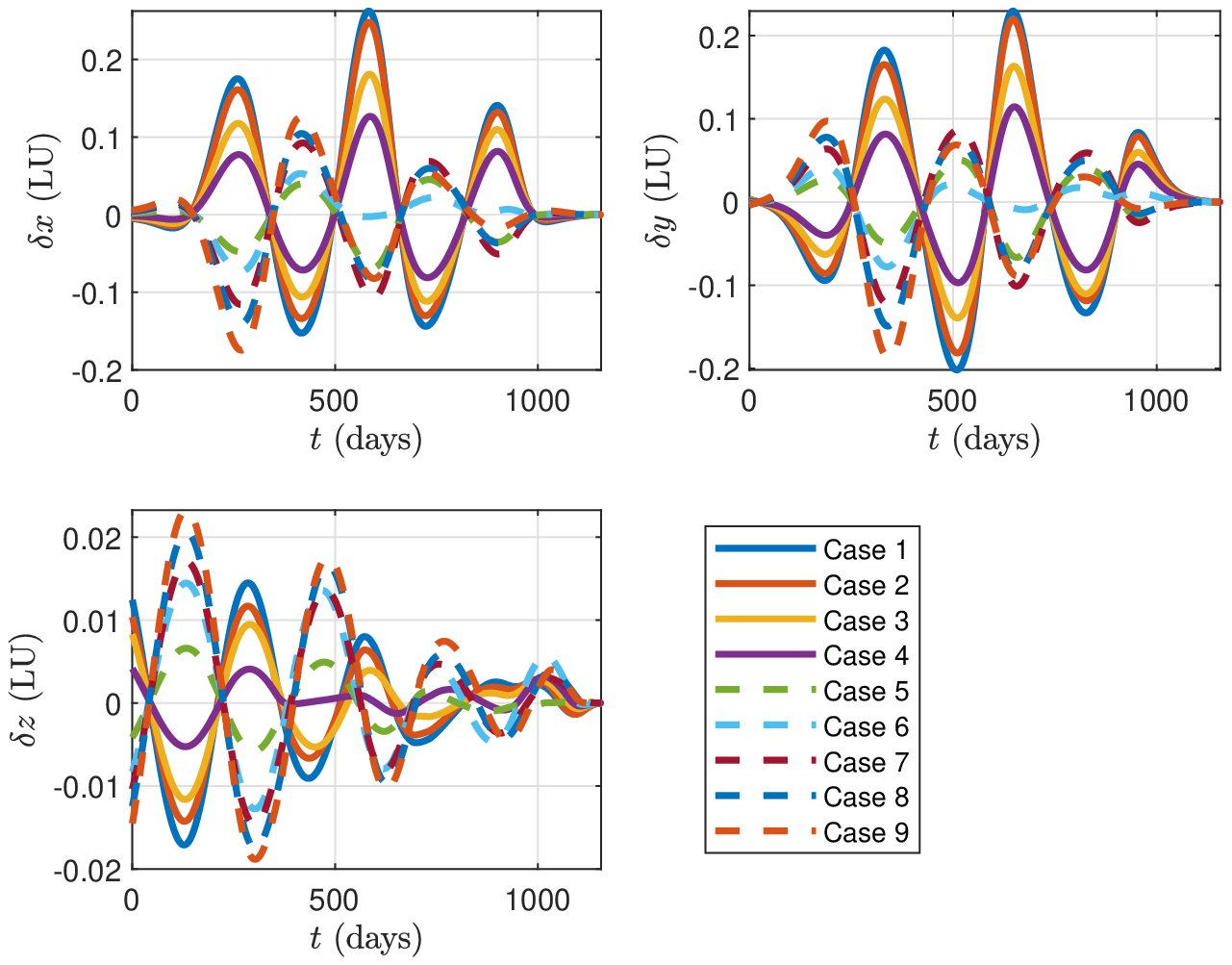}}
	\caption{Differences on coordinates between converged MPSP trajectories and the nominal trajectory for cases in Tab.~\ref{table:ini_per}.}
	\label{fig:ini_per_state_9_cases}
\end{figure}

The comparisons between the thrust angles of the converged MPSP solutions and the nominal thrust angles are shown in Fig.~\ref{fig:Ini_per_9_cases}, where the variations of thrust angles remain to be smooth. The variations of $\alpha$ for cases 5-9 is more apparent than that of cases 1-4, while the $\beta$ oscillations for all cases remain in the vicinity of nominal solution. The comparisons between thrust sequences of the converged MPSP solutions and the corresponding fuel-optimal thrust sequences are shown in Fig.~\ref{fig:Ini_per_u_9_cases}. Case 4 requires the minimum iterations. In this case, the outer-loop MPSP continuation process is not triggered. Case 9 requires the maximum iterations, since the thrust sequence is changed dramatically compared with nominal thrust sequence. From case 4 to 1, the optimal solutions gradually emerges new coast segments, but the obtained MPSP solutions remain the nominal thrust sequence. From case 5 to 9, the initial conditions are becoming tighter, and more thrust is required to drive the spacecraft to the target. The MPSP solutions and optimal solutions are shown a similar trend which gradually increases the thrust segments and reduces the coast segments. From Tab.~\ref{table:ini_per}, the largest increase of fuel consumption is the case 1, which is around $9\%$, while the minimum increase of fuel consumption is the case $4$, which is only $0.16\%$. From case 5 to 9, even though the proposed algorithm requires more iterations for tighter initial conditions, the fuel consumption is nearly optimal. The differences between converged MPSP solutions and the nominal solution on coordinates are shown in Fig.~\ref{fig:ini_per_state_9_cases}. It is interesting to see that the the differences are symmetric for opposite direction of initial perturbations.

\begin{table}
	\centering
	\caption{Perturbations on initial conditions.}
	\begin{tabular}{ccccccc}
		\toprule 
		\toprule
		Case & $\kappa$ & $\|\vect x_f - \vect x(t_f)\|_2 \ (\rm{km})$ & $\|\vect v_f - \vect v(t_f)\|_2 \ (\rm{km/s})$ & $|\lambda_m(t_f)|$ & Newton's iteration & Fuel Increase ($\%$)\\
		\midrule
		1 & $-3$   & $216.45$ & $4.20 \times 10^{-5}$ & $3.74 \times 10^{-8}$  & $40$ & $9.01$\\
		2 & $-2.5$ & $28.77$  & $5.66 \times 10^{-6}$ & $3.11 \times 10^{-9}$  & $21$ & $6.82$\\
		3 & $-2$   & $495.99$ & $9.04 \times 10^{-5}$ & $3.26 \times 10^{-7}$  & $6$  & $2.64$\\
		4 & $-1$   & $328.34$ & $6.07 \times 10^{-5}$ & $1.28 \times 10^{-7}$  & $5$  & $2.27$\\
		5 & $1$    & $128.41$ & $2.43 \times 10^{-5}$ & $2.48 \times 10^{-8}$  & $40$ & $0.16$\\
		6 & $2$    & $9.24$   & $1.86 \times 10^{-6}$ & $1.09 \times 10^{-9}$  & $48$ & $5.02$\\
		7 & $2.5$  & $13.18$  & $2.55 \times 10^{-6}$ & $4.74 \times 10^{-10}$ & $67$ & $1.11$\\
		8 & $3$    & $40.11$  & $1.66 \times 10^{-5}$ & $1.82 \times 10^{-7}$  & $98$ & $2.54$\\
		9 & $3.5$  & $454.25$ & $8.13 \times 10^{-5}$ & $2.41 \times 10^{-7}$  & $135$ &$1.94$\\
		\bottomrule
		\bottomrule
	\end{tabular}
	
	\label{table:ini_per}
\end{table}

\subsection{Perturbations on Terminal Conditions}

Different variations on the terminal positions are simulated to test the proposed method. The totally 8 perturbed terminal positions lie in the vertex of the cube where the reference position locates in the center of the cube, and the side of the cube is set to be $0.04\ {\rm LU}$. The corresponding $8$ cases are shown in Tab.~\ref{table:ter_per_cases}. The simulation results including the terminal errors, the Newton's iteration in the inner-loop MPSP algorithm and the percentage of fuel increase w.r.t the corresponding optimal solutions are shown in Tab.~\ref{table:ter_per}. From the results, it can be observed that all obtained terminal values are within the tolerance.

\begin{figure}
	\centering
	\subfloat[]{\includegraphics[width=0.5\textwidth]{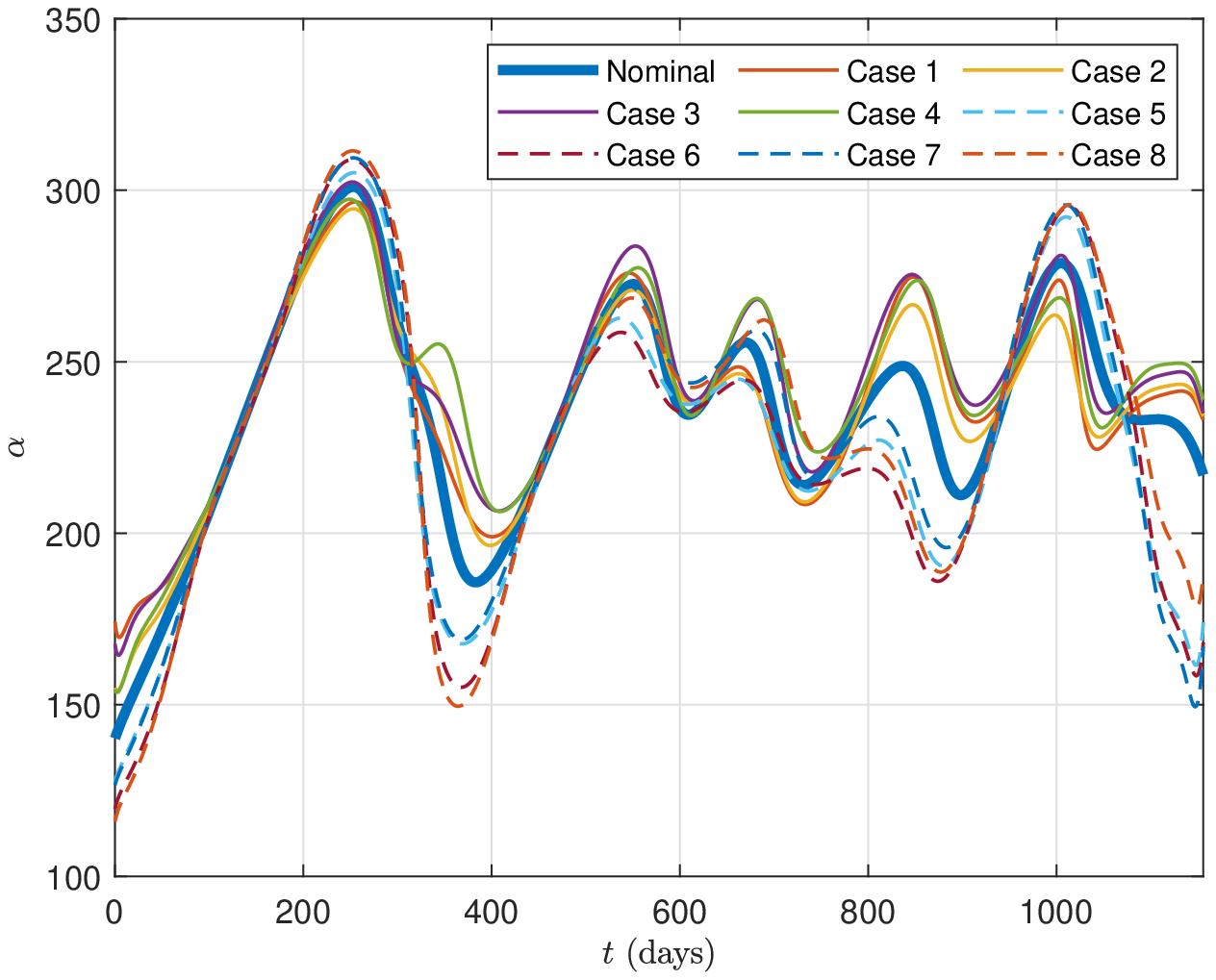}\label{fig:ter_per_alpha_8_cases}}
	\hfill
	\subfloat[]{\includegraphics[width=0.5\textwidth]{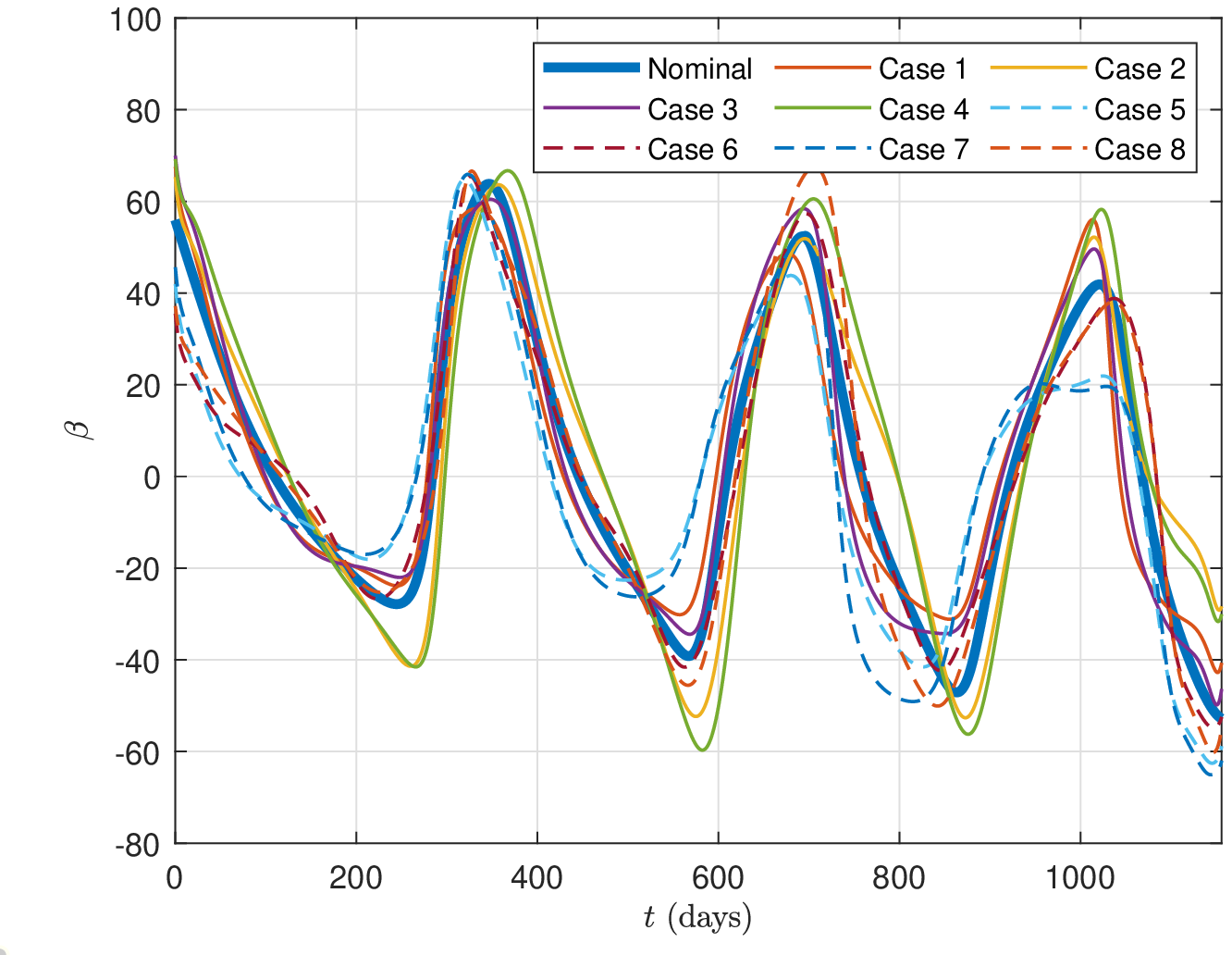}\label{fig:ter_per_beta_8_cases}}
	
	\caption[Simulation results]{The in-plane angle $\alpha$ and out-of-plane angle $\beta$ of the converged MPSP trajectories for variations of terminal positions in Tab.~\ref{table:ter_per_cases}.}
	\label{fig:ter_per_8_cases}
\end{figure}

\begin{figure}
	\centering
	{\includegraphics[width=0.7\textwidth]{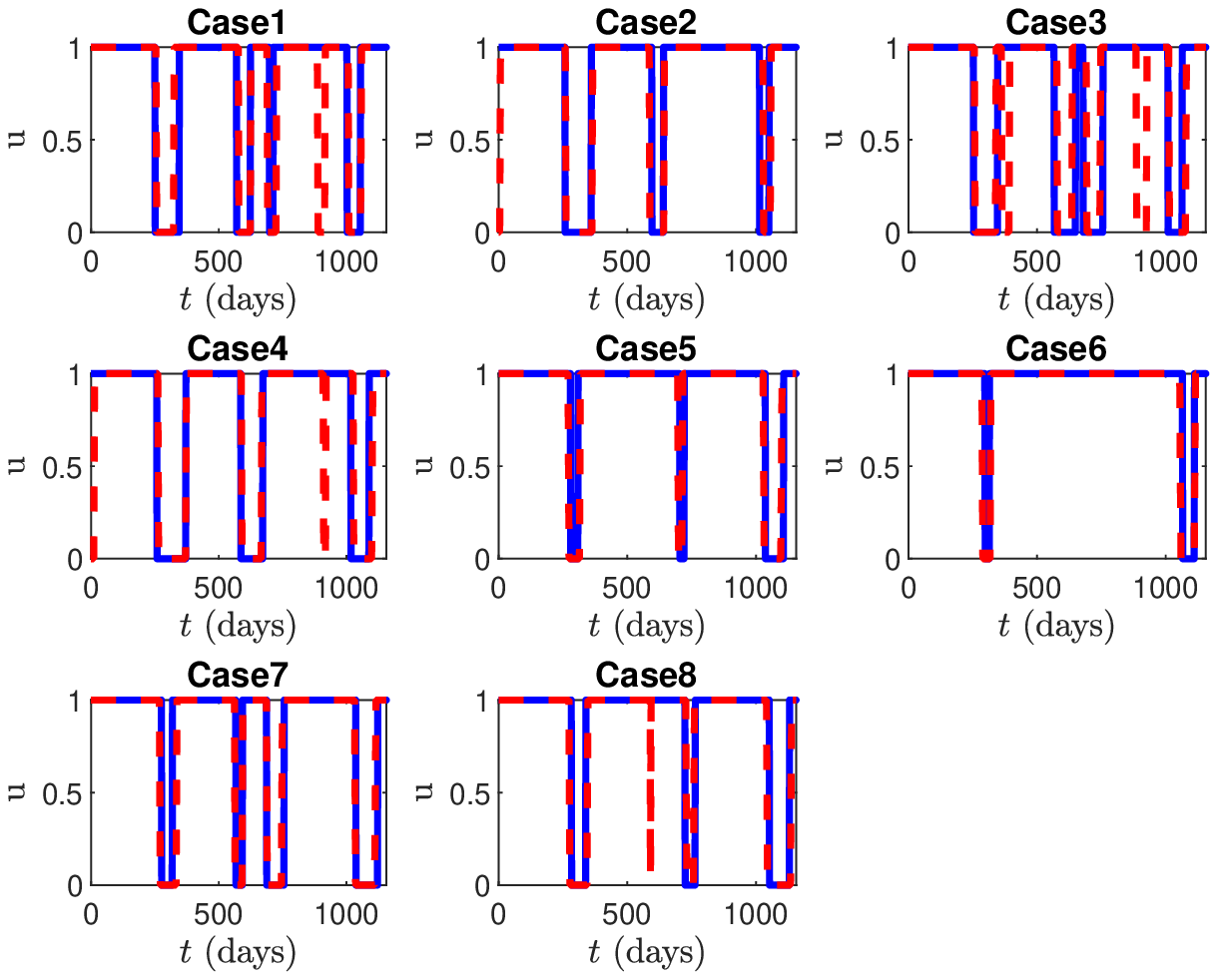}}
	\caption{The comparison of optimal thrust sequences and thrust sequences of MPSP solutions for cases in Tab.~\ref{table:ter_per_cases}. The red dash line: optimal thrust sequences; Blue line: thrust sequences of MPSP solutions.}
	\label{fig:ter_per_u_8_cases}
\end{figure}

\begin{figure}
	\centering
	{\includegraphics[width=0.7\textwidth]{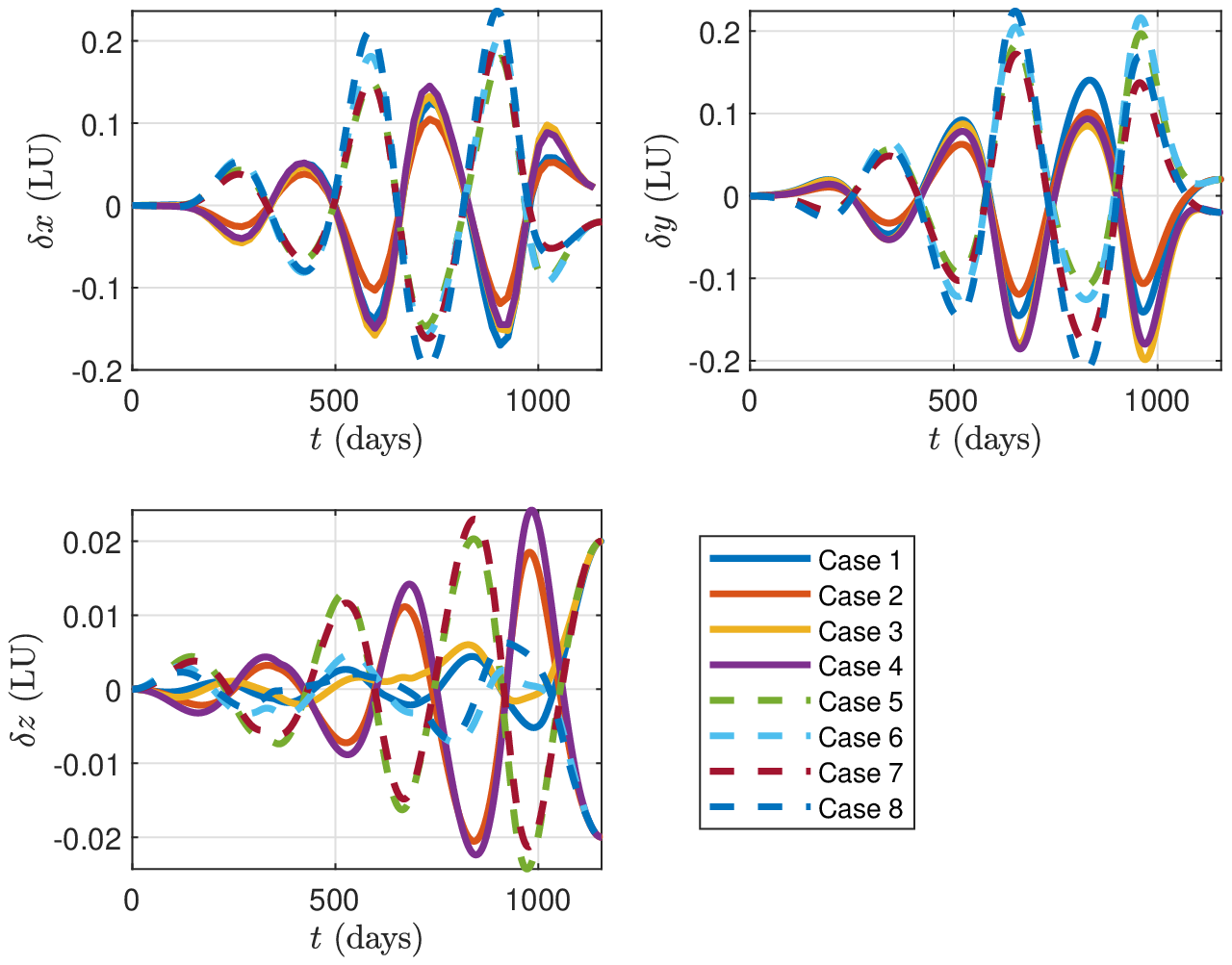}}
	\caption{Differences on coordinates between converged MPSP trajectories and the nominal trajectory for terminal position variation cases in Tab.~\ref{table:ter_per_cases}.}
	\label{fig:ter_per_xyz_8_cases}
\end{figure}

Fig.~\ref{fig:ter_per_8_cases} represents the variations of thrust angles $\alpha$ and $\beta$ w.r.t time. It can be seen that the obtained angle variations are close to the nominal case. Since the terminal $\vect x_f$ for cases 1-4 are farther than cases 5-8 w.r.t the Sun, the $\alpha$ starts more inward in cases 1-4 than cases 5-8. The similar trend is seen from $\beta$ variation. The comparison between the thrust sequences for the converged MPSP solutions and the corresponding optimal thrust sequences are shown in Fig.~\ref{fig:ter_per_u_8_cases}. It is nice to see that the converged MPSP solutions coincide with the optimal solutions well in most cases. For cases 1 and 3, optimal thrust sequences have more coast segments than nominal one. MPSP solutions in these cases prefer to maintain to be the same as the nominal thrust sequence. For cases 2, 4 and 8, MPSP solutions capture the main structure of optimal thrust sequences except some near-impulse thrust segments. In cases 5, 6 and 7, MPSP solutions perfectly coincide with the optimal thrust sequences. In cases 1, 3 and 7 where MPSP solutions remain the nominal thrust sequence, correspond to just 5-6 Newton's iterations. In cases 2, 4, 5 and 8, MPSP requires around 30 Newton's iterations since one less thrust segment is required with $N_{seg,tol} = 2$. In case 6, MPSP requires more iterations because two less coast segment is required with $N_{seg,tol} = 4$. However, for all cases, the fuel consumption is very close to the optimal solution. The minimum fuel consumption is the case 2, which is only $0.79\%$ more than its optimal fuel consumption. The maximum fuel consumption is case 6, which is only $2.36\%$ more than its optimal fuel consumption. The fuel consumption remains nearly optimal even though the thrust sequence is changed w.r.t the nominal thrust sequence. The difference between nominal solution and the converged MPSP solutions are shown in Fig.~\ref{fig:ter_per_xyz_8_cases}. The $x$ difference for case 1-4 and case 5-8 are near symmetric. The $y$ differences are also shown similar symmetry except the last $200$ days. The $z$ differences are shown more complexity and the magnitude of the differences tends to amplify.

\begin{table}
	\centering
	\caption{Cases study for Perturbations on terminal positions.}
	\begin{tabular}{cccc}
		\toprule
		\toprule
		Case & $\delta x_f$ (LU) & $\delta y_f$ (LU) & $\delta z_f$ (LU)\\
		\midrule
		1 & $0.02$  & $0.02$  & $0.02$ \\
		2 & $0.02$  & $0.02$  & $-0.02$\\
		3 & $0.02$  & $-0.02$ & $0.02$\\
		4 & $0.02$  & $-0.02$ & $-0.02$\\
		5 & $-0.02$ & $0.02$  & $0.02$ \\
		6 & $-0.02$ & $0.02$  & $-0.02$\\
		7 & $-0.02$ & $-0.02$ & $0.02$ \\
		8 & $-0.02$ & $-0.02$ & $-0.02$\\
		\bottomrule
		\bottomrule
	\end{tabular}
	
	\label{table:ter_per_cases}
\end{table}

\begin{table}
  	\centering
  	\caption{Simulation results for perturbations on terminal positions.}
  	\begin{tabular}{cccccc}
  		\toprule
  		\toprule
  		Case & $\|\vect x_f - \vect x(t_f)\|_2 \ (\rm{km})$ & $\|\vect v_f - \vect v(t_f)\|_2 \ (\rm{km/s})$ & $|\lambda_m(t_f)|$ & Newton's iteration & Fuel Increase($\%$)\\
  		\midrule
  		1 & $7.98$   & $1.13 \times 10^{-6}$ & $8.35\times 10^{-8}$ & $6$  & $1.33$\\
  		2 & $12.71$  & $2.41 \times 10^{-6}$ & $2.49\times 10^{-9}$ & $40$ & $0.79$\\
  		3 & $8.57$   & $1.60 \times 10^{-6}$ & $1.14\times 10^{-9}$ & $6$  & $2.26$\\
  		4 & $18.85$  & $3.10 \times 10^{-6}$ & $2.52\times 10^{-8}$ & $45$ & $1.01$\\
  		5 & $288.30$ & $5.05 \times 10^{-5}$ & $1.31\times 10^{-7}$ & $32$ & $1.26$\\
  		6 & $16.28$  & $2.82 \times 10^{-6}$ & $2.13\times 10^{-9}$ &$71$  & $2.36$\\
  		7 & $55.82$  & $1.00 \times 10^{-5}$ & $1.26\times 10^{-9}$ &$5$   & $1.28$\\
  		8 & $53.00$  & $8.75 \times 10^{-6}$ & $1.57\times 10^{-9}$ &$40$  & $2.08$\\
  		\bottomrule
  		\bottomrule
  	\end{tabular}
  	
  	\label{table:ter_per}
\end{table}


\subsection{Perturbations on Thruster Parameters}

The perturbations on the thruster parameters are simulated. Specifically, the perturbations on the $T_{\rm max}$ are tested. The percentage of the perturbation w.r.t the nominal solution is set to be $\rm{\eta} = [-10\%,\ -6\%,\ -3\%,\ 3\%,\ 6\%,\ 10\%]$, which corresponds to 6 simulation cases as shown in Tab.~\ref{table:para}. It is assumed that the percentage of the error of $T_{\rm max}$ remains to be the same throughout the flight. The simulation results are summarized in Tab.~\ref{table:para} which gives the terminal errors, the Newton's iteration in the inner-loop MPSP algorithm and the percentage of fuel increase w.r.t the corresponding optimal solutions. For all cases, MPSP algorithm converges successfully.

\begin{figure}
	\centering
	\subfloat[]{\includegraphics[width=0.5\textwidth]{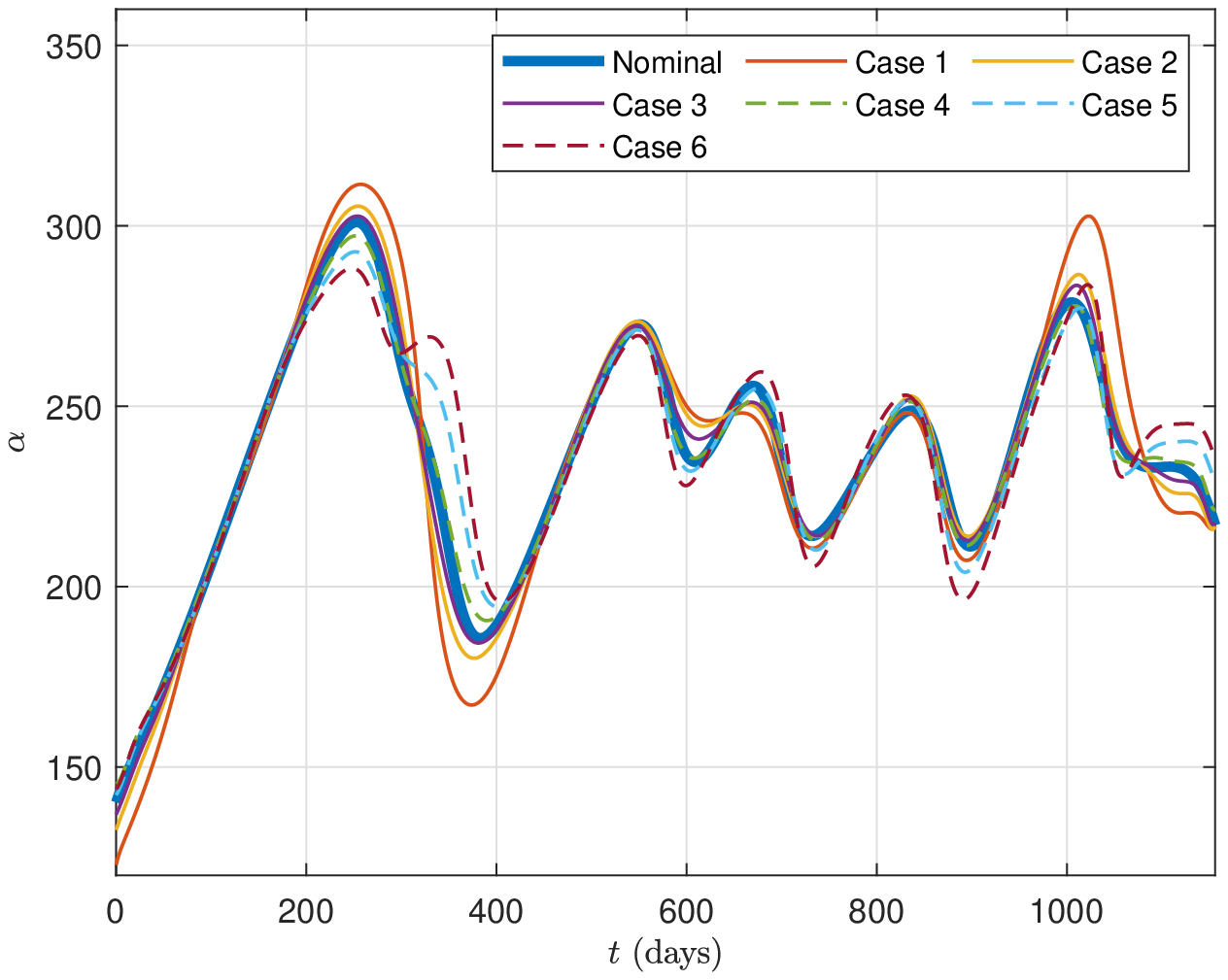}\label{fig:para_per_alpha_6_cases}}
	\hfill
	\subfloat[]{\includegraphics[width=0.5\textwidth]{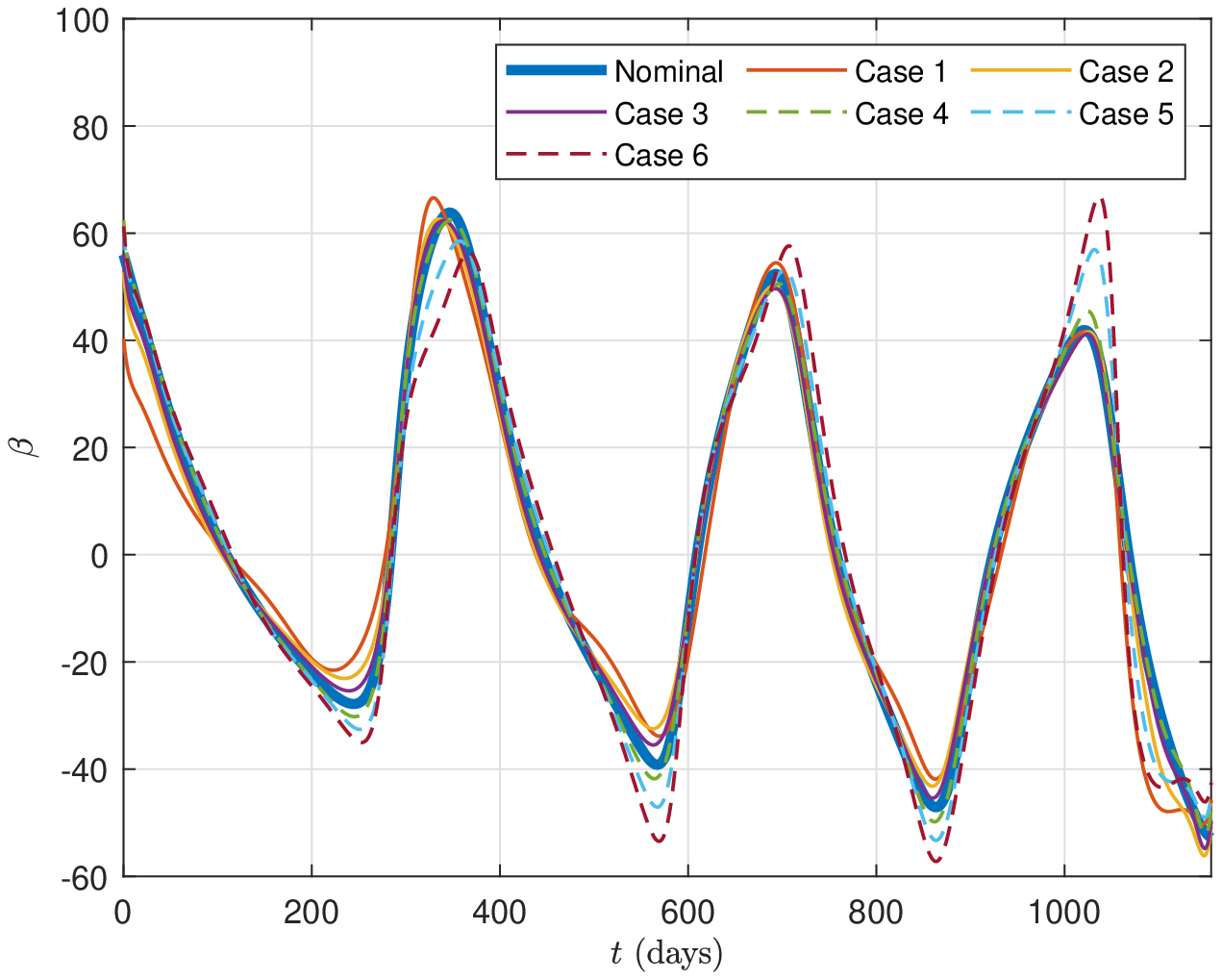}\label{fig:para_per_beta_6_cases}}
	
	\caption[Simulation results]{The variations of the in-plane angle $\alpha$ and out-of-plane angle $\beta$ for the converged MPSP trajectories for variations of $T_{\rm max}$ in Tab.~\ref{table:para}.}
	\label{fig:para_per_6_cases}
\end{figure}

\begin{figure}
	\centering
	{\includegraphics[width=0.7\textwidth]{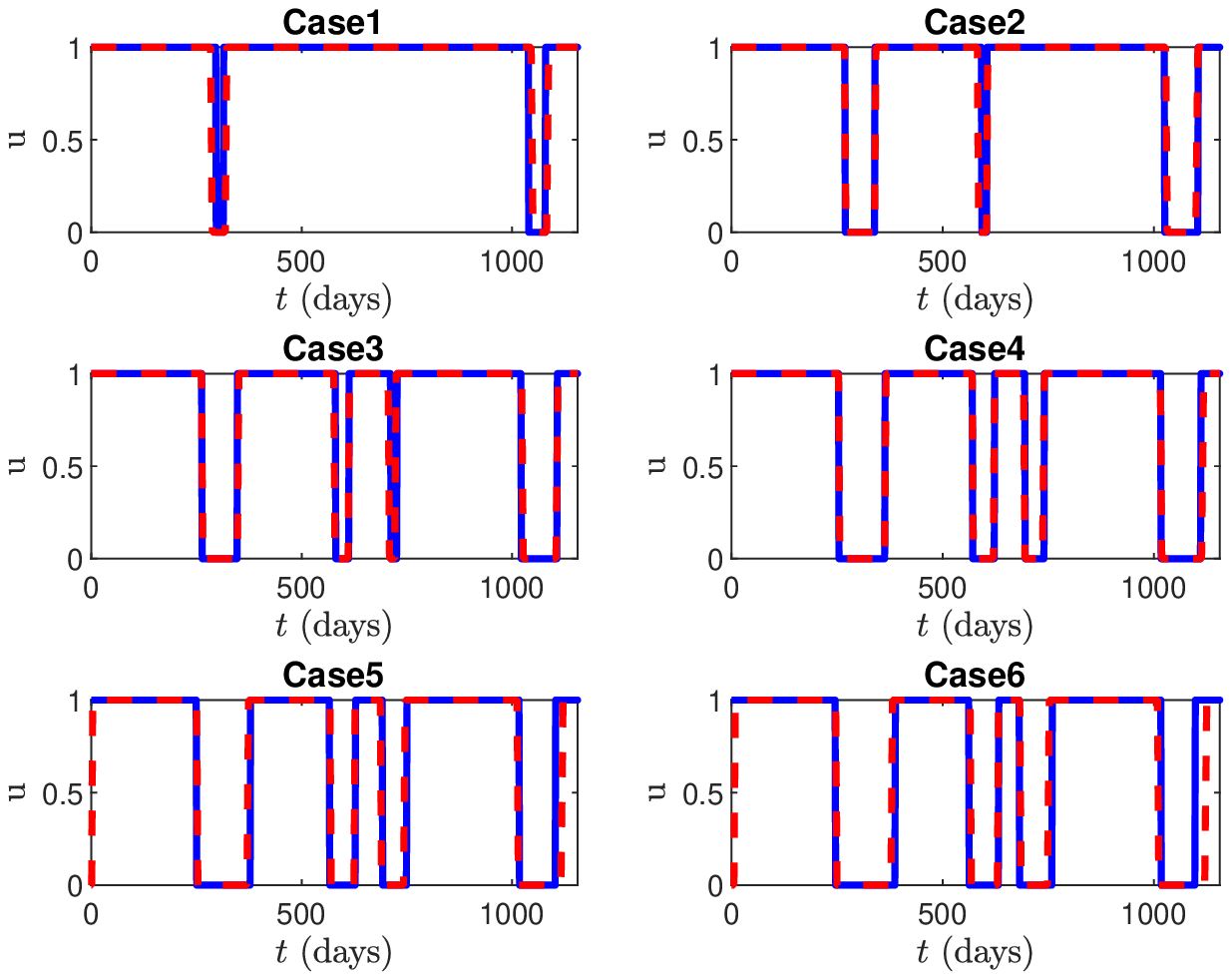}}
	\caption{The comparison of optimal thrust sequences and thrust sequences of MPSP solutions for cases in Tab.~\ref{table:para}. The red dash line: optimal thrust sequences; Blue line: thrust sequences of MPSP solutions.}
	\label{fig:para_per_u_6_cases}
\end{figure}

\begin{figure}
	\centering
	{\includegraphics[width=0.7\textwidth]{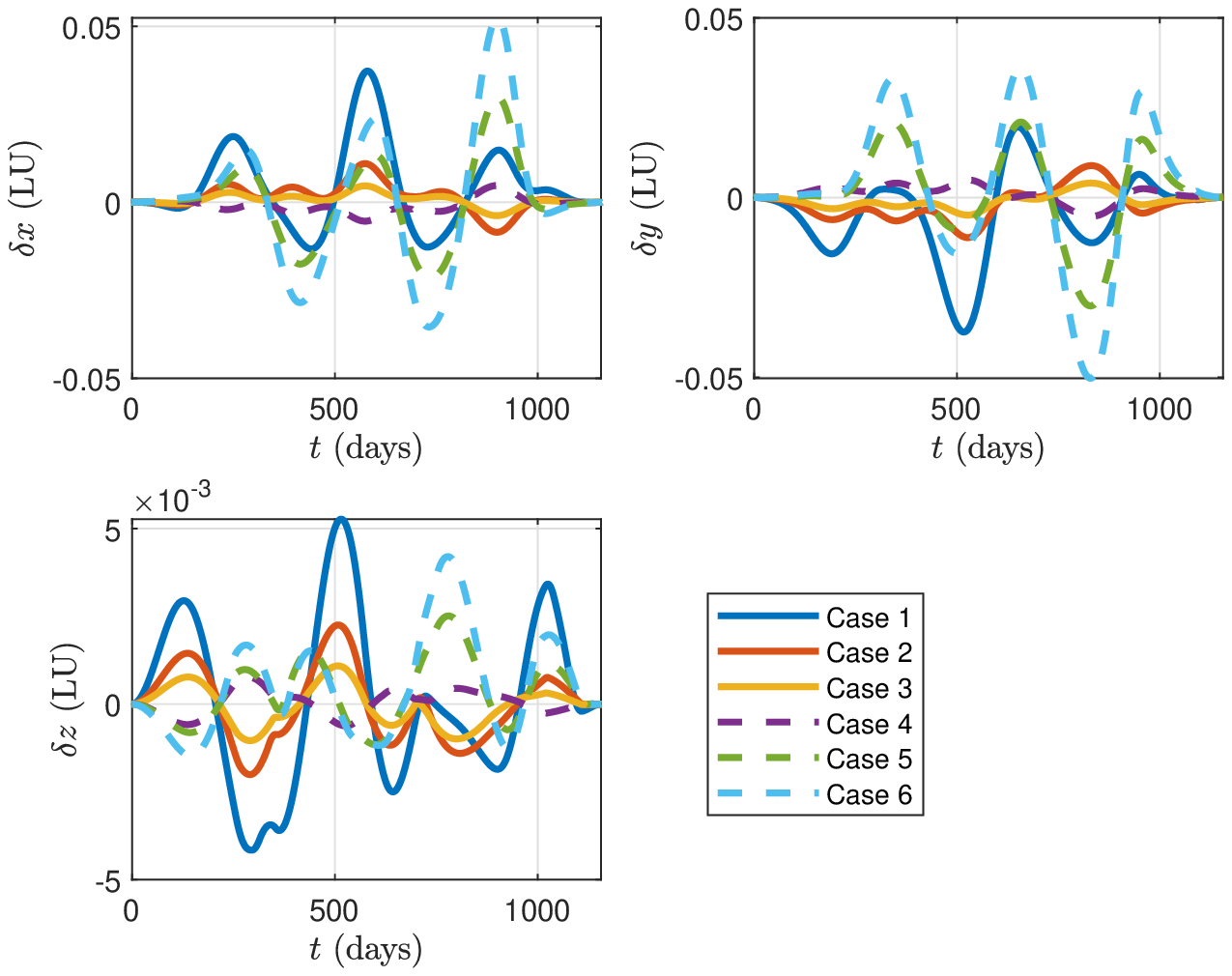}}
	\caption{Differences on coordinates between converged MPSP trajectories and the nominal trajectory for cases of perturbations on $T_{\rm max}$ in Tab.~\ref{table:para}.}
	\label{fig:para_per_xyz_6_cases}
\end{figure}

In Fig.~\ref{fig:para_per_6_cases}, it is nice to see from that the variations of $\alpha$ and $\beta$ are in the vicinity of nominal values. Fig.~\ref{fig:para_per_u_6_cases} illustrates the comparison between the optimal thrust sequences and thrust sequences of the MPSP solutions. The thrust sequences coincide well even when the thrust sequence is changed. The increase of the fuel consumption w.r.t the corresponding fuel-optimal solutions are negligible. As expected, the maximum Newton's iteration occurs in case 1 since the variation of the thrust sequence is the largest. It is also noticed that it only requires 4 to 5 iterations when the thrust sequence remains to be the same as the nominal one. Fig.~\ref{fig:para_per_xyz_6_cases} depicts the trajectory differences between converged MPSP trajectories w.r.t the nominal trajectory on coordinates. The differences for case 1 and 6 are the most obvious since the perturbations on $T_{\rm max}$ are the largest. Not like Figs.~\ref{fig:ini_per_state_9_cases} and \ref{fig:ter_per_xyz_8_cases}, the differences are not symmetry for the cases such as case 1 and 6.

The outcome of this simulation study indicates: 1) the thrust angles of the converged MPSP trajectories remain smooth; 2) the thrust sequence of the MPSP solution prefers to remain to be the nominal thrust sequence when coast segments can be increased. 3) the thrust sequence of MPSP solution can capture the main structure of the optimal thrust sequence when the coast sequence is reduced; 4) Even though the fuel consumption is not included inside the performance index, the MPSP trajectories are competitive in terms of fuel consumption even when the thrust sequence is changed. 

\begin{table}
	\centering
	\caption{Perturbations on $T_{\rm max}$.}
	\begin{tabular}{ccccccc}
		\toprule
		\toprule
		Case & $\eta \ (\%)$ & $\|\vect x_f - \vect x(t_f)\|_2 \ (\rm{km})$ & $\|\vect v_f - \vect v(t_f)\|_2 \ (\rm{km/s})$ & $|\lambda_m(t_f)|$ & Newton's iteration & Fuel Increase ($\%$)\\
		\midrule
		1 & $-10$ & $1.77$   & $3.57 \times 10^{-7}$ & $2.42\times 10^{-9}$ & $93$ & $1.27$\\
		2 & $-6$  & $373.68$ & $5.21 \times 10^{-5}$ & $5.62\times 10^{-7}$ & $24$ & $0.18$\\
		3 & $-3$  & $267.52$ & $5.18 \times 10^{-5}$ & $2.76\times 10^{-7}$ & $4$  & $0.082$\\
		4 & $3$   & $6.50$   & $1.10 \times 10^{-6}$ & $3.15\times 10^{-9}$ & $4$  & $0.062$\\
		5 & $6$   & $12.25$  & $2.10\times 10^{-6}$  & $5.06\times 10^{-9}$ & $5$  & $0.48$\\
		6 & $10$  & $55.35$  & $9.47\times 10^{-6}$  & $1.27\times 10^{-8}$ & $5$  & $1.19$\\
		\bottomrule
		\bottomrule
	\end{tabular}
	
	\label{table:para}
\end{table}

\section{Conclusion}

Unlike the applications with continuous control profile, this paper shows that the sensitive matrix is discontinuous at the bang-bang switching point. A robust two-loop MPSP algorithm is further designed as the low-thrust guidance scheme. Numerical simulations illustrate that the proposed MPSP algorithm is robust for various kinds of perturbations. Besides, the fuel consumption is near optimal even when the thrust sequence is required to be changed. The future work will refine the algorithm design to reduce the total iterations.

\section*{Acknowledgment}
Yang Wang acknowledges the support of this work by the China Scholarship Council (Grant no.201706290024).
\section*{Appendix}
\subsection{Trajectory Discretization}
Instead of discretization using Euler method, a higher order method is used to increase the accuracy. The classical 4th order Runge--Kutta formula is used as
\begin{equation}
\label{eq:Runge-Kutta}
    \left\{
    \begin{array}{ll}
        \vect x_{n+1} &= \vect x_n + \dfrac{h}{6}\left(\vect K_1 + 2\vect K_2 + 2\vect K_3 +\vect K_4\right) \\
        \vect K_1 &= \mathcal{F} (t_{n,1},\vect x_{n,1}), \quad t_{n,1} = t_n, \ \vect x_{n,1} = \vect x_n\\
        \vect K_2 &= \mathcal{F} (t_{n,2},\vect x_{n,2}), \quad t_{n,2} = t_n + \dfrac{h}{2}, \ \vect x_{n,2} = \vect x_n + \dfrac{h}{2}\vect K_1\\ 
        \vect K_3 &= \mathcal{F} (t_{n,3},\vect x_{n,3}), \quad t_{n,3} = t_n + \dfrac{h}{2}, \ \vect x_{n,3} = \vect x_n + \dfrac{h}{2}\vect K_2\\
        \vect K_4 &= \mathcal{F} (t_{n,4},\vect x_{n,4}), \quad t_{n,4} = t_n + h, \ \vect x_{n,4} = \vect x_n + h\vect K_3
    \end{array}
    \right.
\end{equation}

Combining Eq.\ \eqref{eq:Runge-Kutta} with Eq.\ \eqref{eq:expr_x_k_plus_1} yields
\begin{equation}
    \vect F_n(t_n,\vect x_n) = \vect x_n + \dfrac{h}{6}\left(\vect K_1 + 2\vect K_2 + 2\vect K_3 +\vect K_4\right)
\end{equation}
and its differential w.r.t $\vect x_n$ is
\begin{equation}
    \dfrac{\di \vect F_n(t_n,\vect x_n)}{\di \vect x_n} = \vect I_n + \dfrac{h}{6}\left( \dfrac{\di \vect K_1}{\di \vect x_n} + 2\dfrac{\di \vect K_2}{\di \vect x_n} + 2\dfrac{\di \vect K_3}{\di \vect x_n} + \dfrac{\di \vect K_4}{\di \vect x_n} \right)
\end{equation}
where
\begin{equation}
    \dfrac{\di \vect K_i}{\di \vect x_n} = \dfrac{\di \mathcal{F} (t_{n,i},\vect x_{n,i})}{\di \vect x_{n,i}} \dfrac{\di \vect x_{n,i}}{\di \vect x_n},\quad i = 1,2,3,4
\end{equation}

\bibliography{Reference}

\end{document}